\documentclass[12pt]{amsart}
\usepackage[utf8]{inputenc}
\usepackage[T1]{fontenc}
\usepackage[all]{xy}
\usepackage{amscd}
\usepackage{amssymb}
\usepackage{enumerate}
\usepackage{amsthm}
\theoremstyle{plain}
\newtheorem{thm}{Theorem}[section]
\newtheorem{prop}[thm]{Proposition}

\newtheorem{cor}[thm]{Corollary}

\newtheorem{lem}[thm]{Lemma}
\newtheorem{lemdefn}[thm]{Lemma + definition}

\newtheorem{rem}[thm]{Remark}
\newtheorem{defn}[thm]{Definition}
\theoremstyle{definition}

\newtheorem{para}[subsubsection]{}

\newcommand{\Cal}{\mathcal}

\input cyracc.def \font\tencyr=wncyr10 \def\russe{\tencyr\cyracc} \def\Sha{\text{\russe{Sh}}}

\define\cC{\Cal C}

\def\I{\cal I}

\def\Ltens{{\mathop\otimes\limits^{L}}}
\def\L{\Lambda}

\newcommand{\Gm}{{\mathbb{G}_m}}

\newcommand{\N}{{\mathbb{N}}}
\newcommand{\Q}{{\mathbb{Q}}}
\newcommand{\Z}{{\mathbb{Z}}}
\newcommand{\Fp}{{\mathbb F_p}}

\newcommand{\Zp}{{\mathbb Z_p}}

\newcommand{\upc}[1]{\overset {\lower 0.3ex \hbox{${\;}_{\circ}$}}{#1}}

\newcommand{\dlim}{\varinjlim_n}
\newcommand{\ilim}{\varprojlim_n}

\newcommand{\La}{{\operatorname{\Lambda}}}

\newcommand{\cal}{\mathcal}

\renewcommand{\theenumi}{\roman{enumi}}
\newcommand{\limi}{\displaystyle{\lim_{\longrightarrow}}}
\newcommand{\limp}{\displaystyle{\lim_{\longleftarrow}}}

\def\diagram#1{\def\normalbaselines{\baselineskip=0pt\lineskip=10pt\lineskiplimit=1pt}
  \begin{matrix}#1\end{matrix}}
\def\hfl#1#2{\smash{\mathop{\hbox to
10mm{\rightarrowfill}}\limits^{\scriptstyle#1}_{\scriptstyle#2}}}
\def\hflrev#1#2{\smash{\mathop{\hbox to
10mm{\leftarrowfill}}\limits^{\scriptstyle#1}_{\scriptstyle#2}}}
\def\hflcourte#1#2{\smash{\mathop{\hbox to
3mm{\rightarrowfill}}\limits^{\scriptstyle#1}_{\scriptstyle#2}}}
\def\hflrevcourte#1#2{\smash{\mathop{\hbox to
3mm{\leftarrowfill}}\limits^{\scriptstyle#1}_{\scriptstyle#2}}}
\def\injfl#1#2{\smash{\mathop{\hookrightarrow}\limits^{\scriptstyle#1}_{\scriptstyle#2}}}
\def\surjfl#1#2{\smash{\mathop{\twoheadrightarrow}\limits^{\scriptstyle#1}_{\scriptstyle#2}}}
\def\vfl#1#2{\llap{$\scriptstyle #1$}\left\downarrow\vbox to
6mm{}\right.\rlap{$\scriptstyle #2$}}
\def\vflup#1#2{\llap{$\scriptstyle #1$}\left\uparrow\vbox to
6mm{}\right.\rlap{$\scriptstyle #2$}}
\def\vflcourte#1#2{\llap{$\scriptstyle #1$}\left\downarrow\vbox to
2mm{}\right.\rlap{$\scriptstyle #2$}}
\def\vflupcourte#1#2{\llap{$\scriptstyle #1$}\left\uparrow\vbox to
2mm{}\right.\rlap{$\scriptstyle #2$}}

\def\diagram#1{\def\normalbaselines{\baselineskip=0pt\lineskip=10pt\lineskiplimit=1pt}
\begin{matrix}#1\end{matrix}}

\def\debrom{
\makeatletter
\renewcommand{\theenumi}{(\roman{enumi})}
\renewcommand{\labelenumi}{\theenumi}
\makeatother\begin{enumerate}}
\def\finrom{\end{enumerate}}

\def\G{\Gamma}
\def\Zp{{\mathbb Z_p}}
\def\Qp{{\mathbb Q_p}}
\def\mod{{\hbox{-}mod}}
\def\L{\Lambda}

\def\petit{\vspace{.2cm}}

\def\grand{\vspace{.6cm}}

\def\s{\sharp}
\def\O{\cal O}

\begin{document}

\title[] {On the non commutative Iwasawa Main Conjecture for abelian varieties over function fields.}

\author[Trihan]{Fabien Trihan}
\address{Department of Information and Communication Sciences\\
Sophia University\\
Chiyoda-ku, Tokyo, 102-0081\\
Japan}
\email{f-trihan-52m@sophia.ac.jp}

\author[Vauclair]{David Vauclair}
\address{Departement de Mathematiques\\
Universite de Caen\\
Esplanade de la Paix / BP 5186 / 14032 Caen Cedex France}
\email{vauclair@math.unicaen.fr}



\begin{abstract} We establish the Iwasawa main conjecture for semistable abelian varieties over a function field of characteristic $p$ under certain restrictive assumptions. Namely we consider $p$-torsion free $p$-adic Lie extensions of the base field which contain the constant $\Zp$-extension and are everywhere unramified.
Under the usual $\mu=0$ hypothesis, we give a proof which mainly
relies on the interpretation of the Selmer complex in terms of $p$-adic cohomology
\cite{TV} together with the trace formulas of \cite{EL1}.
\end{abstract}
\maketitle
\tableofcontents

\newpage

\section{Introduction}

\subsection{Statement of the main theorem} ~~  \\

Consider a function field $K$ of characteristic $p$, $C$ the corresponding proper smooth geometrically irreducible curve over a finite field $k$ and an abelian variety $A/K$ with N\'eron model $\cal A/C$. Assume for simplicity that $A$ has good reduction everywhere and that the Hasse-Weil $L$ function of $A/K$ does not vanish at $s=1$. In this situation the BSD conjecture predicts that each group involved in the right hand term below is finite and that the following formula holds (where $\hat A$ denotes the dual abelian variety):  \begin{eqnarray}\label{BSD1} L(A/K,1)={\# \Sha(A/K)\over \# A(K)\# \hat A(K)}. {\# H^0(C,Lie(\cal A))\over \# H^1(C,Lie(\cal A))}\end{eqnarray}
Let us restrict our attention to the $p$-adic valuations on both sides.
There are natural perfect complexes of $\Zp$-modules $R\Gamma(C,T_p(\cal A))$, $R\Gamma(C,Lie(\cal A))$ whose cohomology groups are the $p$-part of the groups appearing above. These complexes thus become acyclic after extension of scalars to $\Qp$ hence produce classes in $K_0(\Zp,\Qp)\simeq \Z$. Looking at $L(A/K,1)$ as an element of $K_1(\Qp)\simeq \Qp^\times$ the $p$-part of (\ref{BSD1}) simply becomes \begin{eqnarray} \label{BSD2} \partial (L(A/K,1))=[R\Gamma(C,T_p(\cal A))]+[R\Gamma(C,Lie(\cal A))]\end{eqnarray} where $\partial:K_1(\Qp)\to K_0(\Zp,\Qp)$ is the connecting map of lower $K$-theory (ie. the $p$-adic valuation $v_p:\Qp^\times \to \Z$).
As explained in \cite{KT} the complexes appearing on the right hand side can be related to $p$-adic  (crystalline and then rigid) cohomology and the trace formulae which are known in this context are then sufficient to actually prove (\ref{BSD2}). \\

The main purpose of this paper is to establish a similar statement in the setting of non commutative Iwasawa theory. We don't assume  that $L(A/K,1)\ne 0$ anymore and $A$ is now allowed to have semistable reduction at some given set $Z$ of points of $C$. Consider a Galois extension $K_\infty/K$ which contains the constant $\Zp$-extension $Kk_\infty/K$, is unramified everywhere, and whose Galois group $G=\ilim G_n$ is $p$-adic Lie without $p$-torsion. In that situation, a general result explained in section 2 will allow us to form perfect complexes of modules over the Iwasawa algebra $\L(G):=\ilim \Zp[G_n]$:  $$\begin{array}{rcl} N_{K_\infty}&:=&R\ilim R\Gamma^{Z_n}(C_n,T_p(\cal A)) \\ L_{K_\infty}&:=&R\ilim R\Gamma^{Z_n}(C_n,Lie(\cal A))\end{array}$$
where $R\G^{Z_n}(C_n,-)$ is the functor of cohomology vanishing at $Z_n$ (see. section 4 for a precise definition and the relation to usual Selmer complexes) and is designed to take out the contribution of $Z$. A significant difference with the BSD statement is that here, the cohomology of these complexes are expected to be torsion $\Lambda(G)$-modules even if $L(A/K,1)=0$. This is well known if $K=Kk_\infty$ (see. Cor. \ref{torsionar} or \cite{LLTT} Cor. 2.1.5). 
To go from this case to the general one, we follow the strategy of \cite{CFKSV}, which unfortunately requires an extra assumption, namely the usual $\mu=0$ hypothesis. More precisely, we will prove the following generalization of  \cite{LLTT}, where only the case $K=Kk_\infty$ was considered.

\begin{thm} \label{thmintro} (Thm. \ref{main-theorem}) Let $A/K$, $K_\infty/K$ and $G$ as above. If the $\mu$-invariant
of the Pontrjagin dual of the discrete Selmer group of $A$ over $Kk_\infty$ is trivial, then \debrom 
\item The cohomology $\Lambda(G)$-modules of $N_{K_\infty}$ and $L_{K_\infty}$ are $S^*$-torsion, where $S^*$ is the canonical Ore set of \cite{CFKSV}. The latter complexes thus produce classes $[N_{K_\infty}]$ and $[L_{K_\infty}]$ in $K_0(\L(G),\L(G)_{S^*})$ (denoted  $K_0(\mathfrak{M}_H(G))$ in \emph{loc. cit.}). 
. 
\item There exists a canonical element ${\cal L}_{A/K_\infty}\in K_1(\Lambda(G)_{S^*})$ satisfying \begin{eqnarray}\label{IMCintro1}\partial {\cal L}_{A/K_\infty}=[N_{K_\infty}]+[L_{K_\infty}]\end{eqnarray} where $\partial:K_1(\Lambda(G)_{S^*})\to K_0(\L(G),\L(G)_{S^*})$ denotes the connecting map in lower $K$-theory.
\item The element ${\cal L}_{A/K_\infty}$ verifies the interpolation property \begin{eqnarray}\label{IMCintro2} \rho(\cal L_{A/K_\infty})=L(U,A/K,\rho^\vee,1)\end{eqnarray}
for each Artin $\O$-valued representations $\rho$ of $G$ ($\O/\Z_p$, a totally ramified extension), where 
$\rho(-):K_1(\Lambda(G)_{S^*})\to L^\times\cup \{0,\infty\}$ denotes the corresponding evaluation map, $(-)^\vee$ denotes the contragredient, and $L(U,A/K,\rho,s)$ is the $\rho$-twisted Hasse-Weil $L$-function of $A/K$ without Euler factors at $Z=C - U$.
\finrom 
\end{thm}
 
We also give a similar statement involving complexes $N_{K_\infty}^*$ and $L_{K_\infty}^*$ which are $\La(G)$-dual to $N_{K_\infty}$ and $L_{K_\infty}$ (see Prop.  \ref{main1} $(ii)$, as well as Prop.  \ref{La-duality}, (ii) for the precise statement). \\


Here, the construction of $\cal L_{A/K_\infty}$, the proof that $N_{K_\infty}$ is torsion and the proof of (\ref{IMCintro1}) are simultaneous and rely essentially on the main result of \cite{TV} (a sheafified version of \cite{KT} Prop. 5.13) which yields distinguished triangles of perfect complexes of $\Lambda(G)$-modules relating flat and crystalline cohomology (see Sect. \ref{paradtinfty} and Rem. \ref{remTV}).
The proof of (\ref{IMCintro2}) relies on the comparison of crystalline and rigid cohomology \cite{LST} together with the trace formula for the latter and the codescent properties of the Iwasawa complexes along the tower $K_\infty/K$ (Prop. \ref{La-duality}, relying on Thm. \ref{main1}). 
\\

\subsection{Outline of the paper} ~~ \\

Section 2. Defining the complexes $N_{K_\infty}$, $L_{K_\infty}$ or $N_{K_\infty}^*$, $L_{K_\infty}^*$  occurring in the main conjecture (\ref{IMCintro1}) involves forming projective limits along the Galois tower formed by the $C_n$'s or alternatively inductive limits and taking duals. In order to perform these operations, a convenient framework is given by the derived category of normic systems  $D^b({_{\underline G}}(\Zp\hbox{-}mod))$ defined and studied in \cite{Va}. A normic system is a collection $M_n$ of $G_n$-modules together with equivariant morphisms $M_n\to M_{m}$, $M_m\to M_n$ satisfying a natural compatibility (Def. \ref{defnormic}). The purpose of this section is to show that the collection of derived functors $R\G(C_n,-)$ comes from a functor $R\G(\underline C,-)$ with values in $D^b({_{\underline G}}(\Zp\hbox{-}mod))$.

We place ourself in a setting which is general enough to handle the various cohomology theories (flat, \'etale, crystalline...) involved here. Namely we show (Lem.-Def. \ref{defGnorm}) that in a ringed topos $(E,R)$ any pro-torsor $\underline X=(X_n,G_n)_n$ gives rise to a functor \begin{eqnarray}\label{normicsections}\Gamma(\underline X,-):Mod(E,R)\to {_{\underline G}(\Gamma(E,R)\hbox{-}mod)}\end{eqnarray}  from the category of $R$-modules of $E$ to that of normic systems of $\Gamma(E,R)$-modules.
Using \cite{Va}, the final subsection collects the basic properties (descent, codescent, perfectness, duality, under suitable assumptions, see. Thm. \ref{main1}) of the Iwasawa complexes defined using the functor $$D^+(E,\Zp)\to D(\Lambda(G))$$ obtained by composing $R\Gamma(\underline X,-)$ with $R\ilim$, or alternatively with $\dlim$ and then $RHom(-,\Zp)$. \\

Section 3. We review the basic facts from $K$-theory which are needed for our purpose. This involves mainly the determinant functor for perfect complexes from \cite{Kn} and its behaviour with respect to distinguished triangles and localization following \cite{FK}. We also discuss the evaluation map which appears in (\ref{IMCintro2}).  \\



Section 4. We begin by recalling the definition of the derived functor $R\G^{Z}(C,.)$ of global sections vanishing at a closed subscheme $Z$ which naturally appears in the comparison theorem of \cite{TV} and is denoted $R\Gamma_{ar,V}$ in \cite{KT} (see Rem. \ref{SelGar} for a more precise statement). Next we compare the complex of normic systems $R\G^{\underline Z}(\underline C,T_p(\cal A))$ underlying $N_{K_\infty}$ to the usual Selmer complex of $A/K$. We take the opportunity to give a tractable definition for normic Selmer complexes and prove the expected duality theorem in this setting.

Finally we recall the comparison result of \cite{TV} (which takes place in the small \'etale topos of $C$) and write down the 
fundamental distinguished triangles that follow from it, using (\ref{normicsections}): \begin{eqnarray}\label{tdfundnorm} {\small\begin{array}{c}
R\G(\underline C^\sharp/\Zp,Fil^1D_{log}(\cal A)(-Z))\mathop{\to}\limits^{\bf 1} R\G(\underline C^\sharp/\Zp,D_{log}(\cal A)(-Z))\to R\G(\underline C,Lie(\cal A)(-Z))\mathop{\to}\limits^{+1}\\
R\G^{\underline Z}(\underline C,T_p(\cal A))\to R\G(\underline C^\sharp/\Zp,Fil^1D_{log}(\cal A)(-Z))\mathop{\to}\limits^{\bf 1-\phi} R\G(\underline C^\sharp/\Zp,D_{log}(\cal A)(-Z))\mathop{\to}\limits^{+1}\end{array}}\end{eqnarray}
Here $\phi$ is a semi-linear map such that $\phi{\bf1}=p Frob$. By applying $R\ilim$, these in turn yield distinguished triangles of perfect complexes of $\Lambda(G)$-modules satisfying the derived codescent property. This will be the main ingredient for the proof of \ref{thmintro}.  \\

Section 5. We put everything together in order to prove the Iwasawa main conjecture. The third term of the first distinguished triangle above  is a $k$-vector space and the arrow denoted $\bf 1$ thus becomes invertible after inverting $S^*$. Whence an endomorphism $(\bf 1-\phi)_{S^*}\bf 1_{S^*}^{-1}$ acting on the localization $P_{K_\infty,S^*}$ of $$P_{K_\infty}:=R\limp R\G(C_n^\sharp/\Zp,D_{log}(\cal A)(-Z))$$
In the case where $G=\G$, the base change formula in crystalline cohomology together with a semi-linear argument shows that the first term in the second distinguished triangle above vanishes as well after $S^*$-localization. In the general case, an argument using Nakayama's lemma ensures that it is still the case under the $\mu=0$ assumption. The endomorphism $({\bf 1}-\phi)_{S^*}{\bf 1}_{S^*}^{-1}$ is thus invertible, allowing us to  \emph{define} $\cal L_{A/K_\infty}\in K_1(\Lambda_{S^*})$ as its determinant, rendering (\ref{IMCintro1}) almost tautological. Let us hint the idea of the proof of (\ref{IMCintro2}). Using the descent properties of the normic section functor together with the comparison between crystalline and rigid cohomology from \cite{LST}, we prove an isomorphism \begin{eqnarray}\label{515}L\Ltens_{\La_O(G)}(
V_{\rho, \La_{O}(G)})^*\Ltens_{\La(G)}P_{K_\infty}\simeq R\G_{rig,c}(U/L,pr^*D^\dagger\otimes U(\rho^\vee))\end{eqnarray} where $U(\rho^\vee)$ denotes the unipotent convergent isocrystal associated to $\rho$ (Prop. \ref{descrig}). On the one hand the trace formula in rigid cohomology shows that $L(U, A/K,\rho^\vee,1)$ coincides with the alternated product of the determinants of $1-p^{-1}Frob$ acting on the $L$-valued cohomology  of the right hand side. Since we have no morphism from $\Lambda_{S^*}$ to $L$ we may not use directly the base change property of the $Det$ functor together with a localized version of (\ref{515}) in order to relate the action of $1-p^{-1}Frob$ on the right hand side with $\rho(\cal L)$. We turn this difficulty by investigating carefully the spectral sequence of codescent along $Kk_\infty/K$  with the definition of the evaluation map in mind.

\bigskip
\noindent{\bf Acknowledgements.} The first author is supported by JSPS. Both authors thank the referees for their careful reading and suggestions to simplify the paper.

\section{The normic section functor} \label{SAC} ~~ \\

The purpose of this section is to show that under reasonable conditions the cohomology of a topos along a Galois tower with group $G$ naturally gives rise to perfect complexes of $\L(G)$-modules satisfying natural properties such as derived descent, codescent and duality along the tower. This goal is achieved in Thm. \ref{main1}. \\ ~~ \\

\subsection{Normic systems} \label{NS} ~~ \\

We define and study briefly the category of normic systems along a profinite group $G$. The following definitions slightly generalize those given in \cite{Va}.



\begin{defn} \label{defnormic} Consider a profinite group $G$, $N$ a filtered index set, $(H_n)_{n\in N}$ the filtered set
of the normal open subgroups and denote $G_n=G/H_n$, $G_{m,n}=H_n/H_m$, $n\le m$. Given a commutative ring $S$.

$(i)$ ${_{G}Mod(S)}$ the category of $S$-modules endowed with a discrete action of $G$, ie. which are the union of their fixed points by the $G_n$'s.

$(ii)$ ${_{G_n}Mod(S)}$ the full subcategory of ${_{G}Mod(S)}$ whose objects
are those on which $H_n$ acts trivially. For $m\ge n$, the inclusion functor ${_{G_n}Mod(S)}\rightarrow {_{G_m}Mod(S)}$ has a left and a
right adjoint, described respectively as the coinvariants functor $(-)_{G_{m,n}}$ and the fixed points functor $(-)^{G_{m,n}}$.


$(iv)$ ${_{\underline G}Mod(S)}$ the category of $Mod(S)$-valued normic systems is defined a follows:

- an object is a triple $\underline M=((M_n)_{n\in N}, (j_{n,m}: M_n\rightarrow M_m)_{n\le m}, (k_{m,n}:M_m\rightarrow M_n)_{n\le m})$ satisfying the
following properties:

\item $(Norm1)$ $M_n$ is an object of ${_{G_n}Mod(S)}$, $j_{n,m}$ and
$k_{m,n}$ are morphisms of ${_{G}Mod(S)}$.

\item $(Norm2)$ If $n_1\le n_2\le n_3$ then $j_{n_2,n_3}\circ j_{n_1,n_2}=j_{n_1,n_3}$ and $k_{n_2,n_1}\circ k_{n_3,n_2}=k_{n_3,n_1}$.

\item $(Norm3)$ If $n\le m$, $j_{n,m}\circ k_{m,n}:M_m\rightarrow M_m$ coincides with the endomorphism $\sum_{g\in G_{m,n}}c_g$ of $M_m$.
Note that this is not only a morphism in $Mod(S)$, but also in ${_{G}Mod(S)}$ since $H_n\lhd G$.

- a morphism $\underline \phi:\underline M\rightarrow \underline M'$ is a collection of morphisms of ${_{G}Mod(S)}$, $(\phi_n:M_n\rightarrow M'_n)_{n\in
N}$  such that the following squares commute for each couple $n\le m$:
$$\xymatrix{M_m\ar[r]^{\phi_m}&M_m'&\hspace{2cm}&M_m\ar[d]_{k_{m,n}}\ar[r]^{\phi_m}&M_m'\ar[d]_{k_{m,n}}\\
M_n\ar[u]^{j_{n,m}}\ar[r]^{\phi_n}&M_n'\ar[u]^{j_{n,m}}&\hspace{1.5cm}&M_n\ar[r]^{\phi_n}&M_n'}$$
\end{defn}

Let us make some remarks about the category ${_{G}Mod(S)}$ of discrete $G$-objects.

- For $M,M'\in _{G}Mod(S)$: $$Hom_{_{G}Mod(S)}(M,M')\simeq \limp_n\limi_mHom_{_{G}Mod(S)}(M^{H_n},M'^{H_m}).$$

- The category $_{G}Mod(S)$ is a full subcategory of $Mod(S[G])$. The inclusion functor has a right adjoint $M\mapsto \cup_n M^{H_n}$. In particular, $_{G}Mod(S)$ has enough injectives since $Mod(S[G])$ has.
\medskip


%
%
%
%
%
%
%

We now turn to some properties of the category $_{\underline G}Mod(S)$ of normic systems.

\begin{prop} Small inductive and projective limits exist in $_{\underline G}Mod(S)$ and commute to the component functors $_{\underline G}Mod(S)\to _{\underline G_n}Mod(S)$. In particular, $_{\underline G}Mod(S)$ is an abelian category.
\end{prop}
Proof. In order to form the inductive (resp. projective) limit indexed by a set $I$ in ${_{\underline G}\mathcal C}$, it suffices to form the limit of the components and endow them with the $j_{n,m}$ and $k_{m,n}$ provided by functoriality. This prove the first statement. The second one follows since being abelian is a property of limits (\cite{KS}, Def. 8.2.8 and Def. 8.3.5).  \begin{flushright}$\square$\end{flushright}

The following lemma,  inspired by \cite{TW} has been pointed out
by B. Kahn.

\begin{lem}

$(i)$ If $G$ is in fact a finite group, then the category ${_{\underline G}Mod(S)}$ is equivalent to $Mod(\mu_\lhd(G,S))$ where
$\mu_\lhd(G,S)$ is the \emph{normal Mackey algebra}, defined as the quotient of the free associative algebra $S\{\{c_g,g\in G\},\{j_{n,m},
k_{m,n},n\le m\}\}$ by the following relations:

- $c_{g_1}c_{g_2}=c_{g_1g_2}$, $c_gj_{n,m}=j_{n,m}c_g$, $c_gk_{m,n}=k_{m,n}c_g$, for $g,g_1,g_2\in G$ and $n\le m$ in $N$.

- $j_{n_2,n_3}j_{n_1,n_2}=j_{n_1,n_3}$, $k_{n_2,n_1} k_{n_3,n_2}=k_{n_3,n_1}$, for $n_1\le n_2\le n_3$ in $N$.

- $c_gj_{n,n}=j_{n,n}$,  $j_{n,n}=k_{n,n}$, for $n\in N$ and $g\in H_n$.

- $j_{n,m}k_{m,n}=\sum_{g\in G_{m,n}}c_g j_{m,m}$, for $n\le m$ in $N$.

- $\sum_{n\in N}j_{n,n}=1$.

$(ii)$ In general, ${_{\underline G}Mod(S)}$ is equivalent to the category of cocartesian sections of the cofibered category above $N^{op}$
associated to the covariant pseudo-functor on $N^{op}$ mapping $n$ to $Mod(\mu_\lhd(G_n,S))$, and $m\ge n$ to the functor $Mod(\mu_\lhd(G_m,S))\rightarrow Mod(\mu_\lhd(G_n,S))$ given by extension of scalars through the morphism $\mu_\lhd(G_m,S)\rightarrow \mu_\lhd(G_n,S)$ defined by killing the $j_{a,b}$'s and $k_{b,a}$'s with $b\nleq n$.
\end{lem}

Proof. $(i)$  Let us describe the functor ${_{\underline G}Mod(S)}\rightarrow Mod(\mu_\lhd(G,S))$. A normic system $\underline M$ is
sent to $\oplus_nM_n$, together with its obvious structure of $\mu_\lhd(G,S)$-module: $c_g$ acts on every components whereas $j_{n,m}$ (resp. $k_{m,n}$)
sends the $n$-th (resp. $m$-th) component into the $m$-th (resp. $n$-th) one. By definition, a morphism of normic systems $\underline M\rightarrow
\underline M'$ consists in a collection of morphisms $\phi_n:M_n\rightarrow M_n'$, compatible with the $c_g$'s, the $j_{n,m}$'s and the $k_{m,n}$'s.
It thus gives rise to a morphism $\oplus_n M_n\rightarrow \oplus_nM'_n$, compatible with the action of $\mu_\lhd(G,S)$ and this is clearly compatible
to composition. We thus have defined the desired functor ${_{\underline G}Mod(S)}\rightarrow Mod(\mu_\lhd(G,S))$. It now remains to
notice that both the full faithfulness and the essential surjectivity of this functor immediately follow from the isomorphism of algebras
$\mu_\lhd(G,S)\simeq\prod_{n\in N}j_{n,n}\mu_\lhd(G,S)$.

$(ii)$ We apply $(i)$ to the case $G=G_n$. Letting $n$ vary, this gives an equivalence of cofibered categories. Whence the result, since
${_{\underline G}Mod(S)}$ is clearly equivalent to the category of cocartesian sections of the cofibered category associated to $n\mapsto
{_{\underline{G_n}}Mod(S)}$.
\begin{flushright}$\square$\end{flushright}

The following corollary answers a question of \cite{Va}.
\begin{cor} \label{enoughinj}
The category ${_{\underline G}Mod(S)}$ has enough injectives.

\end{cor}
Proof.  Let $\mathcal F/N^{op}$ denote the fibered category $n\mapsto Mod(\mu_\lhd(G_n,S))$ considered previously and let $Sect(\mathcal F/N^{op})$
(resp. $Cocart(\mathcal F/N^{op})$) denote its category of sections (resp. cocartesian sections). Denoting $|N|$ the discrete category underlying $N$,
there are three obvious forgetful functors
$$\begin{array}{ccccccc}Cocart(\mathcal F/N^{op})&\mathop\rightarrow\limits^{F_1}&Sect(\mathcal F/N^{op})
&\mathop\rightarrow\limits^{F_2} &\prod_{|N|}Mod(\mu_\lhd(G_n,S))&\mathop\rightarrow\limits^{F_3} &Mod(S)^{|N|}\end{array}$$ each of
which is exact and whose composition is faithful. The result will follow formally if one proves that each of them has a right adjoint. For $F_3$, this is easy and left to
the reader.
For $F_2$, this results from \cite{SGA4}, Vbis, 1.2.10 (note that here $\mathcal F/N^{op}$ is indeed bifibered). For $F_1$, we use the following general easy
fact concerning an arbitrary cofibered category $\mathcal F/\mathcal S$ with cocleavage $f\mapsto f_!$:

\emph{Fact}: Assume that for any $f:X\rightarrow Y$ the following properties are verified:

- The functor "composition with $f$" :$\mathcal S_{/X}\rightarrow \mathcal S_{/Y}$ is cofinal.

- $f_!$ commutes to projective limits indexed by $\mathcal S_{/X}$.

Then the inclusion functor $Cocart(\mathcal F/\mathcal S)\rightarrow Sect(\mathcal F/\mathcal S)$ has a right adjoint, which takes a section
$X\mapsto s(X)$ to the cocartesian section $X\mapsto \limp_{f:Z\rightarrow X} f_!s(Z)$.

\begin{flushright}$\square$\end{flushright}

\begin{rem} One may show that ${_{\underline G}Mod(S)}$ has enough projectives as well.
\end{rem}

\ \\




\subsection{From sheaves to normic systems} ~~ \\

As before let $G=\ilim G_n$ denote a profinite group and let $G_n=G/H_n$. Recall that the classifying topos $B_G$ is the category of discrete left $G$-sets. Consider another arbitrary topos $E$ together with its structural morphism $\pi:E\to Set$ and assume given a projective system $\underline X=(X_n)$ of $G$-objects of $E$ such that $X_n$ is a torsor of $E$ under $\pi^{-1}G_n$.
In this situation we have for $m\ge n$ a commutative diagram of topoi as follows: $$\xymatrix{&&&E\ar[llld]_{p_\infty}\ar[lld]^{p_m}\ar[ld]^{p_n}\ar[d]^\pi\\ B_G\ar[r]_{q_m}&B_{G_m}\ar[r]_{q_{m,n}}&B_{G_n}\ar[r]&Set}$$
Here the horizontal arrows are by functoriality of the classifying topos with respect to the group (their inverse and direct images functors are inflation and fixed points by the adequate subgroup) and the oblique morphisms, defined by $\underline X$, are as follows:

- If $Y$ is in $B_{G_n}$ and $Y^0$ denotes the underlying set, the formula $g.(x,y)=(gx,gy)$ defines an action of $\pi^{-1}G_n$ on $X_n\times \pi^{-1}Y^0$ and $p_n^{-1}Y=\pi^{-1}G_n\backslash (X_n\times \pi^{-1}Y^0)$ is the coinvariant object of this action. If $Z$ is an object of $E$, then $p_{n,*}Z$ is the set $Hom(X_n,Z)$ endowed with the left action of $G_n$ induced by the inverse action on $X_n$: $(g.f)(x)=f(g^{-1}x)$.

- If $Y$ is in $B_G$ then the formula $g.(x,y)=(gx,gy)$ defines an action of $p_\infty^{-1}G$ on  $X_m\times
\pi^{-1}(Y^{H_n})^0$  for each $m,n$ in $N$ and $p_\infty^{-1}Y=\limi_n\limp_m {\pi^{-1}G}\backslash (X_m\times \pi^{-1}(Y^{H_n})^0)$ (Note that $p_\infty^{-1}$ is exact, since
$\limp_m$ is essentially constant while $\limi_n$ is filtered). If $Z$ is an object of $E$ then $p_{\infty,*}Z$ is the set $\dlim Hom(X_n,Z)$ endowed with the left action of $G$ induced by the inverse action of $G$ on $f$.

Endowing the set $G_n$ with its left action by translations turns it into an object of $B_G$ (resp. $B_{G_m}$, if $m\ge n$) whose image by $p_\infty$ (resp. $p_m$) is nothing but $X_n$.


For the purpose of what follows we will consider a ring $S$. Endowing it with a trivial action of $G$, we may as well view it as a ring of $B_{G_n}$ or $B_G$.

\begin{defn} \label{defGnormG} We define a functor
$$\G(\underline G,-):Mod(B_G,S)\to {_{\underline G}Mod(S)}$$
by sending an $S$-module $M$ of $B_G$  to $(M_n,j_{n,m},k_{m,n})$, where $M_n=M^{H_n}$, $j_{n,m}$ is the  inclusion and $k_{m,n}$ the trace map.
\end{defn}

\begin{lemdefn} \label{defGnorm} Let $G$, $S$, $E$ and $\underline X$ be as above. Consider furthermore a ring $R$ of $E$ and a ring homomorphism $S\to \pi_*R$. Letting $Mod(E,R)$ denote the category of $R$-modules of $E$, we define the functor of \emph{Normic sections along $\underline X$} $$\G(\underline X,-):Mod(E,R)\rightarrow {_{\underline G}Mod(S)}$$
by composing $p_{\infty,*}:Mod(E,R)\to Mod(B_G,S)$ with the functor $\G(\underline G,-):Mod(B_G,S)\to {}{_{\underline G}Mod(S)}$ of Def. \ref{defGnormG}. 

For any $F$ in $Mod(E,R)$, the object $((M_n),(j_{n,m}),(k_{m,n})):=\G(\underline X,F)$ satisfies the following properties:

- $M_n$ is the restriction of scalars to $S$ of the $R(X_n)$-module $F(X_n)$ endowed with the action of $G_n$ coming from the right action of $G_n$ on $X_n$: $(g,x)\mapsto g^{-1}x$.

-$j_{n,m}:M_n\to M_m$ is the restriction along $X_m\rightarrow X_n$.
and induces an isomorphism $M_n\simeq M_m^{G_{m,n}}$ (ie. $F(X_n)\mathop\rightarrow\limits_{} F(X_m)^{G_{m,n}}$). 
\end{lemdefn}

Proof. By definition we have $M_n=M^{H_n}=q_{n,*}M$, where $M=p_{\infty,*}F$. Now, $q_np_\infty\simeq p_n$, therefore $M_n\simeq p_{n,*} F=F(X_n)$. The second property is clear.
\begin{flushright}$\square$\end{flushright}

\begin{rem} If the topos $E$ is locally connected, it is possible to build trace maps along finite locally free morphisms, such as $X_m\to X_n$. One may then show that $k_{m,n}$ coincides with the trace map. 
\end{rem}

\begin{lem} Consider another topos $E'$ with structural morphism $\pi':E'\to Set$ and a morphism of ringed topoi $f:(E',R')\to (E,R)$. Let us moreover denote $\underline X'=(f^{-1}X_n)$ the projective system of torsors of $E'$ deduced from $\underline X$ and consider the morphism $S\to \pi'_*R'$ induced by $h:S\to \pi_*R$ and $R\to f_*R'$. The functor $\G(\underline X',-)$ associated to these data is subject to a canonical isomorphism $$R\G(\underline X',-)\simeq R\G(\underline X,Rf_*(-))$$ of functors $D^+(E',R')\to D^+({_{\underline G}Mod(S)})$.
\end{lem} 
Proof. This simply follows from the fact that the morphism of ringed topoi $(E',R')\to (B_G,S)$ defined by $\underline X$ coincides with the one obtained by composing $f$ with the morphism $(E,R)\to (B_G,S)$ defined by $\underline X$.
\begin{flushright}$\square$\end{flushright}

The derived functor $$R\G(\underline X,-):D^+(Mod(E,R))\rightarrow D^+({_{\underline
G}Mod(S)})$$ 
will be our tool to deduce the fundamental distinguished triangles (which are the main ingredient in our proof of the main conjecture) from the comparison isomorphism of \cite{TV}. It might be worth to emphasize that the $n$th component functor $$(-)_n:D^+({_{\underline G}Mod(S)})\rightarrow D^+({_{G_n}Mod(S)})$$  sends $R\G(\underline X,F)$ to $R\G(X_n,F_{|X_n})$.


\subsection{Descent, codescent and duality} \label{DCAD}



 \label{Iwcomp} ~~ \\

We now assume that the ring $S$ of Def. \ref{defGnormG} is $\Zp$ and we discuss how to pass from normic systems to Iwasawa modules.


In order to be able to invoke directly the results of \cite{Va}, we now assume that $\La:=\Zp[[G]]:=\limp \Zp[G_n]$ is Noetherian of finite global dimension $d+1$ and has an open pro-$p$-subgroup. As is well known, these
condition is in particular verified if $G$ is a compact $p$-adic Lie group of dimension $d$ without $p$-torsion.
 

\smallskip

- \emph{Limits.}  Forgetting the $k_{m,n}$'s (resp. the $j_{n,m}$'s), inflating from $G_n$-modules to discrete $G$-modules (resp. abstract
$\La$-modules) and then forming the limit of the resulting inductive (resp. projective) system gives rise to a
functor which will abusively be denoted $$\limi:{_{\underline G}Mod(\Zp)}\to {_GMod(\Zp)} \\ \hbox{ (resp. } \limp:{_{\underline G}Mod(\Zp)}\to Mod(\L) \hbox{)}$$

The first is exact and thus passes to derived categories while the second is only left exact, but right derivable. From now on, \emph{we always assume that  $R\limp$ has finite cohomological dimension}. This is e.g. the case if $N$ has a numerable cofinal subset (as is always the case in practice).

\begin{prop}                                                                                                                          \label{adjlim}
There are canonical adjunctions $$\begin{array}{ccccc}(\limi,R\G(\underline G,-))&:&D^+({_GMod(\Zp)})&\rightarrow
&D^+({_{\underline G}Mod(\Zp)})\\
(\underline{\La}\Ltens_{\La}(-),R\limp)&:&D^-({_{\underline G}Mod(\Zp)})&\rightarrow& D^-({\La})\end{array}$$ Here, $\underline{\La}$ denotes the
natural normic system of right ${\La}$-modules $(\Zp[G_n])$ and $\underline{\La}\otimes_{\La}(-):Mod(\La)\rightarrow
{_{\underline G}Mod(\Zp)}$ denotes the induced right exact, left derivable functor. 
\end{prop}
Proof. The derived version are easily deduced from the obvious ordinary adjunctions using that ${_GMod(\Zp)}$ and ${_{\underline G}Mod(\Zp)}$ have enough
injectives (cf. Prop. \ref{enoughinj}). In \cite{Va}, this was not known and a finiteness assumption was thus needed to avoid deriving $\limp$ (cf. loc. cit. Prop. 4.2. and Rem.
4.3).
\begin{flushright}$\square$\end{flushright}

In the above proposition, it is possible to make the adjunction morphisms involved functorial at the level of the complexes. Also, they can
still be provided a functorial cone, as in loc. cit.


- \emph{Duality.} Consider a $\Zp$-module $I$. If $M$ (resp. $M_n$) is a $\Zp$-module endowed with a left action of $G$ (resp. $G_n$) then $Hom_\Zp(M,I)$ (resp. $Hom_\Zp(M_n,I)$) is endowed a left action of $G$ (resp. $G_n$) as well by reversing the action on $M$ (ie. $(g.f)(m)=f(g^{-1}m)$). If the action of $G$ is discrete the action of $G$ on $Hom_\Zp(M,I)$ extends to an action of $\L$. If now $\underline M=(M_n,j_{n,m},k_{m,n})$ is a normic system then $Hom_\Zp(\underline M,I):=(Hom_\Zp(M_n,I),Hom_\Zp(k_{m,n},I),Hom_\Zp(j_{n,m},I))$ is a normic system as well. Using an injective resolution of $\Zp$ we get duality functors $$\begin{array}{ccccc}RHom_\Zp(-,\Zp)&:&D^b({_GMod(\Zp)})&\rightarrow& D^b(\La)\\
RHom_\Zp(-,\Zp)&:&D^b({_{\underline G}Mod(\Zp)})&\rightarrow& D^b({_{\underline G}Mod(\Zp)})\end{array}$$
It might be useful to point out the following relation to Pontryagin duality: $$RHom_\Zp(M\Ltens \mathbb Q_p/\Zp,\mathbb
Q_p/\Zp)\simeq RHom_\Zp(M,\Zp)$$

If now $M$ is a $\L$-module we view $Hom_\L(M,\L)$ as a left $\L$-module via the right action of $\L$ on itself and the involution $g\mapsto g^{-1}$. Using projective resolutions of $M$ we get a functor
$$\begin{array}{ccccc}RHom_{\La}(-,\La)&:&D^b(\La)&\rightarrow & D^b(\La)\end{array}$$

%
%

\begin{thm} \label{main1}  Let $E$, $X_n$ be as in the previous paragraph.

$(i)$ $(Descent)$   There is a canonical isomorphism 
$$R\G(\underline X,-)\mathop\rightarrow\limits^\sim R\G(\underline G,\limi R\G(\underline
X,-))$$
of functors $D^+(E,\Zp)\rightarrow D^+({_{\underline G}Mod(\Zp)})$.

$(ii)$  Consider a bounded complex $F$ of $\Zp$-modules of $E$ (resp. and assume that there exists $q_0$) such that
$H^q(X_n,F_{|X_n})$ is finitely generated over $\Zp$ (resp. trivial) for any $n\in N$, $q\ge 0$ (resp. $q\ge q_0$).

$(Codescent)$ Functorially in $F$, There is a canonical isomorphism in $D^+({_{\underline G}Mod(\Zp)})$:

$$\underline{\La}\Ltens_{\La}R\limp R\G(\underline X,F) \mathop\rightarrow\limits^\sim
R\G(\underline X,F)$$

$(Perfectness)$ $R\limp R\G(\underline X,F)$ and $RHom_\Zp(\limi R\G(\underline X,F),\Zp)$ are in $D^p(\La)$ (the derived category of perfect complexes).

$(Duality)$ Functorially in $F$, There is a canonical isomorphism in $D^p(\La\mod)$:
$$RHom_\Zp(\limi R\G(\underline X,F),\Zp)\simeq RHom_\La(\limp R\G(\underline X,F),\La)$$
\end{thm}
Proof. $(i)$ We are going to check that the adjunction morphism of $(\limi,R\G(\underline G,-))$ is an isomorphism. Replacing $F$ by an injective
resolution and truncating, one reduces to the case where $F$ is an injective $\Zp$-module of $E$ placed in degree $0$. Let $\underline
M:=\G(\underline X,F)$. Then $\forall n$, $M_n$ is injective in ${_{G_n}Mod(\Zp)}$ (indeed $M_n$ corresponds to $p_{n*}F$) and the descent map induces an isomorphism
$M_n\rightarrow M_m^{G_{m,n}}$. Similarly $\dlim M_n$ is injective and the adjunction morphism of descent occurring in  $D^+({_{\underline
G}Mod(\Zp)})$
is represented at level $n$ by the isomorphism $\G(X_n,F_{|X_n})\rightarrow \G(H_n,\limi_m \G(X_m,F_{|X_m}))$.

$(ii)$ Under the stated assumptions,

- $(Codescent)$ and $(Perfectness)$ follow from $(i)$, by \cite{Va} 4.4 $(ii)$ et 4.6 $(i)$.

- $(Duality)$ follows from $(Codescent)$, by loc. cit. 4.6 $(iii)$.
\begin{flushright}$\square$\end{flushright}

\begin{rem} Properties $(Codescent)$ and $(Perfectness)$ are a crucial tool in this paper. It seems plausible that they hold without assuming that $G$ has finite cohomological dimension but we have not tried to check this.\\
\end{rem}  

\subsection{Examples} \label{sectexamples} ~~ \\

Start with a projective system of surjective \'etale schemes $C_n$ over a base scheme $C$ and assume that the $C_n$'s are endowed with compatible actions of the $G_n$'s such that each $G_n\times C_n\to C_n\times_C C_n$, $(g,x)\mapsto (gx,x)$ is an isomorphism. (In view of our applications, let us notice that if $C$ and $C_n$ are proper smooth curves over a field with respective function fields $K$ and $K_n$ then the latter condition holds as soon as the unramified extension $K_n/K$ is Galois). These data represent a projective system of torsors, say $\underline X$ of the small \'etale topos $E=C_{et}$ and we are thus in the situation of the previous paragraphs.
 
Of course these data also produce a projective system of torsors, say $\underline X'$, of the big, say flat, topos $E'=C_{FL}$. Note that if $\epsilon:C_{FL}\to C_{et}$ denotes the natural morphism, then $\underline X'=\epsilon^{-1}\underline X$.

Consider now the small \'etale crystalline topos $E''$ of $C/\Zp$ and the projection morphism $u:E''=(C/\Zp)_{crys,et}\to C_{et}=E$. Then we define a projective system of torsors of $E''$ by pulling back via $u$: $\underline X'':=u^{-1}\underline X$. 

As explained in the previous paragraph, these data give rise to compatible normic section functors, ie. to an essentially commutative diagram $$\xymatrix{D^+(C_{FL},R')\ar[r]^-{R\epsilon_*}\ar[dr]_-{R\G(\underline X',-)}&D^+(C_{et},R)\ar[d]^-{R\G(\underline X,-)}&D^+((C/\Zp)_{crys,et},R'')\ar[l]_-{Ru_*}\ar[dl]^{R\G(\underline X'',-)}\\ &D^+({_{\underline G}Mod(\Zp)})&}$$
whenever one is given $\Zp$-algebras $R$, $R'$, $R''$ in $E$, $E'$, $E''$ and homomorphisms $R\to \epsilon_*R'$, $R\to \epsilon_* R''$. 

The $n$-component of the left (resp. middle, resp. right) normic section functor computes the cohomology of the localized topos $C_{FL}/X'_n$, (resp. $C_{et}/X_n$, resp. $(C/\Zp)_{crys,et}/X_n''$), ie. of $C_{n,FL}$ (resp. $C_{n,et}$, resp. $(C_n/\Zp)_{crys,et}$). We will thus use the more suggestive notations $R\G(\underline C_{FL},-)$ (resp. $R\G(\underline C_{et},-)$, resp. $R\G_{crys,et}(\underline C/\Zp,-)$).


\section{K-theory for $\L(G)$} \label{KT}


For the convenience of the reader, we review the $K$-groups and determinant functors which are needed for
our purpose.
Unless specified otherwise, $R$ denotes a
(unitary) ring.

\subsection{Review of $K_0$ and $K_1$} \label{K0K1} ~~ \\

The following definitions can be found in \cite{Ba}, ch. VII, IX.

- $K_0(R)$ is the abelian group defined by generators $[P]$, where $P$ is a finitely generated projective $R$-modules, and relations
\begin{enumerate}
\item $[P]=[Q]$ if $P$ is isomorphic to $Q$ as a $R$-module.
\item $[P\oplus Q]=[P]+[Q]$
\end{enumerate}

- $K_1(R)$ is the abelian group defined by generators $[P,f]$, where $P$ is a finitely generated projective $R$-module and $f$ is an automorphism of
$P$ with relations (the group law is denoted multiplicatively):
\begin{enumerate}
\item $[P,f]=[Q,g]$ if $P$ is isomorphic to $Q$ as a $R$-module via an isomorphism which is compatible with $f$ and $g$.
\item $[P\oplus Q,f\oplus g]=[P,f][Q,g]$
\item $[P,fg]=[P,f][P,g]$
\end{enumerate}

Note that if the ring $R$ is Noetherian and regular, then forgetting the word ``projective'' does not change the definitions (cf \cite{Ba}, IX,
Proposition 2.1).

- The functor $Hom_R(-,R)$ realizes an anti-equivalence between finitely generated projective left modules and finitely generated projective right modules. As a
result, one gets an isomorphism $K_i(R)\simeq K_i(R^{op})$, if $R^{op}$ denotes the opposed ring. Both this isomorphism and its inverse will be
denoted $(-)^*$ (eg. $[P,f]^*=[P^*,f^*]$ if $f^*$ denotes the transpose of $f$).

- Morita equivalence. Let us fix a finitely generated projective right $R$-module $V$ and let $V^*=Hom_R(V,R)$. Then $V$ (resp. $V^*$) is naturally
endowed with a natural structure of $((End_R(V),R)$-bimodule (resp. $(R, End_R(V))$-bimodule). Since $V$ is projective and finitely generated, one
has a canonical isomorphism of $(R,R)$-bimodules:  $V^*\otimes_{End_R(V)}V\simeq R$ (resp. of $(End_R(V),End_R(V))$-bimodules: $V\otimes_RV^*\simeq
End_R(V)$.

The functor $V^*\otimes_{End_R(V)}(-):End_R(V)\mod\rightarrow R\mod$ is thus an equivalence of categories, with $V\otimes_R(-)$ as quasi-inverse. In
particular, there is a canonical isomorphism $K_i(End_R(V))\simeq K_i(R)$, $i=0,1$. \\ ~~ \\

\subsection{The determinant functor}  ~~ \\

Let $\cal P(R)$ denote the category of strictly perfect (ie. bounded with projective finitely generated objects) complexes, $K^p(R)$ its homotopy category and $D^p(R)$ its essential image in the derived category (which is naturally equivalent, and identified, to $K^p(R)$). If $D^{p}(R)^{is}$ denotes the subcategory of $D^p(R)$ where morphisms are the isomorphisms of the latter we have by \cite{Kn}, th. 2.3, 2.12 a canonical functor \begin{eqnarray} \label{defDet} Det_R:D^{p}(R)^{is}\simeq K^p(R)^{is}\rightarrow \cal C_R\end{eqnarray} where $\cal C_R$ is the Picard category (ie. a category together with an endobifunctor, referred to as the \emph{product}, endowed with associativity and commutativity isomorphisms satisfying natural compatibilities, unit objects and where every object and every
morphism is invertible; cf \cite{Kn} appendix A) considered in \cite{FK} 1.2:
\begin{itemize}

\item An object of ${\cal C}_R$ is a couple $(P,Q)$ of finitely generated projective $R$-modules.

\item $Mor((P,Q),(P',Q'))$ is empty if $[P]-[Q]\neq [P']-[Q']$ in $K_0(R)$. Else, there exists an $R$-module $M$ such that
$$P\oplus Q'\oplus M\simeq P'\oplus Q\oplus M.$$
We set $I_M:=Isom(P\oplus Q'\oplus M,P'\oplus Q\oplus M)$ and $G_M:=Aut(P'\oplus Q\oplus M)$ and we  define the set of morphisms from $(P,Q)$ to
$(P',Q')$ as $K_1(R)\times^{G_M} I_M$ where the right hand side denotes the quotient of $K_1(R)\times I_M$ by the action of $G_M$ given by
$$(x,y)\mapsto (x\bar{g},g^{-1}y))$$
where $x\in K_1(R)$, $y\in I_M$, $g\in G_M$ and $\bar{g}$ is its image in $K_1(R)$. As seen easily, this set does not depend on $M$, up to a
canonical isomorphism, and this fact can be used to define composition in a natural way.  Note that by definition, one has a canonical identification
$Aut_{\cal C_R}((P,Q))=K_1(R)$ for any object $(P,Q)$.  \\

\item The product is defined as
$$(P,Q).(P',Q'):=(P\oplus P',Q\oplus Q').$$ and admits naturally $(0,0)$ as a unit.
Every object $(P,Q)$ in ${\cal C}_R$ admits  $(P,Q)^{-1}:=(Q,P)$ as a natural inverse.

\end{itemize}

It follows immediately from its construction that the functor $(\ref{defDet})$ is compatible with base change, ie. the diagram $$\diagram{D^p(R')&\hfl{Det_{R'}}{}&\cal C_{R'}\cr
\vflup{R'\Ltens_R(-)}{}&&\vflup{}{R'\otimes_R(-)}\cr D^p(R)&\hfl{Det_R}{}&\cal C_R}$$ is naturally pseudo-commutative for any $R'/R$. Let us quickly review the construction of $Det_R$ following \cite{FK} 1.2.

\begin{itemize}
\item A strictly perfect complex $P$ is sent to $Det_R(P):=(P^+,P^-)$ where
$P^+:=\oplus_{i\in 2\mathbb Z}P^i$ and $P^-:=\oplus_{i\in 1+2\mathbb Z}P^i$.

\item Any exact sequence $0\rightarrow P'\rightarrow P\rightarrow P''\rightarrow 0$ of $\cal P(R)$ induces a canonical isomorphism $Det_R(P)\simeq
Det_R(P'). Det_R(P'')$. Using this, one constructs a canonical \emph{trivialization} $can_P:Det_R(0)\simeq Det_R(P)$, for any strictly perfect
complex $P$ which is acyclic.

\item If $Cone(f)$ denotes the mapping cone of a morphism $f:P\rightarrow P'$ of strictly perfect com\-plexes then  $Hom_{\cal C_R}(Det_R(P),Det_R(P'))$ identifies with $Hom_{\cal
C_R}(Det_R(0),Det_R(Cone(f)))$ (indeed $Cone(f)^+ \simeq P^-\oplus P'^+$ and $Cone(f)^-\simeq P^+\oplus P'^-$). In particular, when
$f$ is a quasi-isomorphism, then $can_{Cone(f)}:Det_R(0)\to Det_R(Cone(f))$ induces a morphism $Det_R(P)\rightarrow Det_R(P')$, which we denote $Det_R(f)$. One may check that it is compatible with composition and only depends on the homotopy class of $f$ so that the functor $Det_R$ is finally defined. \end{itemize}

We will make essential use of the homorphism
\begin{eqnarray} \label{Detaut} Det_R:Aut_{D^p(R)}(P)\to Aut_{C_R}(Det_R(P))=K_1(R)\end{eqnarray}
which is induced by $(\ref{defDet})$ for any perfect complex $P$. Of course if $P$ is reduced to a (finitely generated projective) module $P^0$ placed in degree zero then $Det_R(P)=(P^0,0)$ and $Det_R(f):(P^0,0)\to (P^0,0)$ is the class of $f^0:P^0\to P^0$. In that case $(\ref{Detaut})$ is thus nothing but the tautological map $Aut_{R}(P^0)\to K_1(R)$. \smallskip

Let us now state some multiplicative properties.

\begin{itemize}
\item Consider a morphism between exact sequences of strictly perfect complexes  $$\diagram{0&\hfl{}{}&P_1'&\hfl{}{}&P_1&\hfl{}{}&P_1''&\hfl{}{}&0\cr && \vfl{f'}{}&&\vfl{f}{}&&\vfl{f''}{}\cr
0&\hfl{}{}&P_2'&\hfl{}{}&P_2&\hfl{}{}&P_2''&\hfl{}{}&0}$$ where vertical arrows are quasi-isomorphisms. The
following square commutes: $$\diagram{Det_R(P_1)&\simeq& Det_R(P_1').Det_R(P_1'')\cr\vfl{Det_R(f)}{}&&\vfl{}{Det_R(f').Det_R(f'')}\cr Det_R(P_2)&\simeq&
Det_R(P_2').Det_R(P_2'')}$$

\item Consider an automorphism $f$ of some $P$ in $K^p(R)$. Then:

- $Det_R(f[1])=Det_R(f)^{-1}$ in $Aut_{\cal C_R}(Det_R(P))= K_1(R)$.

- $f$ is the homotopy class of some automorphism of the complex $P$. Let us fix such one and denote $f^q:P^q\to P^q$ its component of degree $q$. Then in $K_1(R)= Aut_{\cal C_R}(Det_R(P))= Aut_{\cal C_R}(Det_R(P^q))$, we have $Det_R(f)=\prod Det_R(f^q)^{(-1)^{q}}$.

- If $P$ is cohomologically perfect (ie. each $H^q$ is an object of $D^p(R)$) then $Det_R(f)=\prod Det_R(H^q(f))^{(-1)^{q}}$ in $K_1(R)= Aut_{\cal C_R}(Det_R(P))= Aut_{\cal C_R}(Det_R(H^q(P)))$.
\end{itemize}


\subsection{Localization}\label{relative K} ~~ \\

The main $K$-group of interest for Iwasawa theory is a relative one. We recall its definition (\cite{FK} 1.3). Consider a strictly full triangulated subcategory $\Sigma$ of $D^p(R)$. The group $K_1(R, \Sigma)$ is then defined by generators and relations as
follows:
\begin{itemize}
\item Generators: $[C, a]$ where $C$ is an object of $\Sigma$ and $a:Det_{R}(0)\rightarrow
Det_{R}(C)$ is a \emph{trivialization of $C$}.

\item Relations: Let $C,C', C"$ be objects of $\Sigma$.
\begin{itemize}
\item If $C\simeq 0$ then
$[C, can_C] = 1$.

\item If $f:C\simeq C'$, then compatible trivializations of $C$ and $C'$ give rise to the
same element in $K_1(R,\Sigma)$ (ie. $[C',Det_R(f)\circ a]=[C,a]$ if $a:Det_R(0)\simeq Det_R(C)$).

\item If $0\rightarrow C'\rightarrow C\rightarrow C''\rightarrow 0$ is an exact sequence of $\cal P(R)$, and
$a:1\simeq Det_{R}(C)$, $a' :1\simeq Det_{R}(C')$, then $$[C,a]=[C',a'].[C'',a'']$$ where $a:1\mathop\rightarrow\limits^{a'.a''}
Det_R(C').Det_R(C'')\simeq Det_R(C)$.
\end{itemize}
\end{itemize}

There is a localization exact sequence (cf \cite{FK} 1.3.15) $$\diagram{K_1(R)&\hfl{}{}& K_1(R,\Sigma)&\hfl{\partial}{}& K_0(\Sigma)\hfl{}{}& K_0(R)}$$
(here $K_0(\Sigma)$ denotes the Grothendieck group of the triangulated category $\Sigma$) where:

- the first map sends $[P,f]\in K_1(R)$ to the complex $[P\mathop\rightarrow\limits^f P]$ placed in degrees $[-1,0]$, together with the trivialization
which is represented by the identity of $P$.

- $\partial$ sends $[C,f]\in K_1(R,\Sigma)$ to the element $-[C]\in K_0(\Sigma)$.

- the last map sends the class of a strictly perfect complex $[C]$ to the alternated sum  $\sum (-1)^i[C^i]$.

As checked easily, this localization sequence is functorial with respect to $(R,\Sigma)$: if $R'$ is an $R$-algebra, and if
$\Sigma'$ is a strictly full triangulated subcategory of $D^p(R')$ containing the essential image of $\Sigma$ under the functor
$R'\otimes^L_R(-):D^p(R)\rightarrow D^p(R')$ then one has a morphism $K_1(R,\Sigma)\rightarrow K_1(R',\Sigma')$ which is compatible with the obvious functoriality maps.
%

\medskip

If now $T\subset R$ is a left denominator set (cf \cite{FK}, 1.3.6) then we can apply the above constructions to the full triangulated subcategory
$\Sigma_T$ of $D^p(R)$ consisting of complexes whose image under the functor $R_T\otimes^L_R(-):D^p(R)\to D^p(R_T)$ become acyclic. We note that if $R$ is Noetherian and regular then $K_0(\Sigma_T)$ is isomorphic to
the Grothendieck group of finitely generated $T$-torsion modules, by sending $[C]$ to $\sum (-1)^iH^i(C)$. For a general $R$ we have the following result.

\begin{prop} (\cite{FK}, 1.3.7) There is a canonical isomorphism of groups $$K_1(R,\Sigma_T)\simeq K_1(R_T)$$
sending $[C,a]$ to the isomorphism $Det_{R_T}(0)\mathop\rightarrow\limits^{a_T}Det_{R_T}C_T\mathop\rightarrow\limits^{can_{C_T}^{-1}}Det_{R_T}(0)$ viewed as an element of $K_1(R_T)=Aut_{\cC_{R_T}}(Det_{R_T}(0))$ (here
$a_T$ and $C_T$ are deduced from $a$ and $C$ by localization). This isomorphism is functorial with respect to $(R,T)$.
\end{prop}
\begin{flushright}$\square$\end{flushright}

Let us mention an alternative characterization of this
isomorphism. Consider an endomorphism $f$ of a strictly perfect complex $P$ such that $f_T$ is a quasi-isomorphism. Identifying $Hom_{\cal C_R}(Det_R(P),Det_R(P))$ with $Hom (Det_R(0),Det_R(Cone(f)))$ sends the identity of $Det_R(P)$ to a morphism $triv: Det_R(0)\rightarrow
Det_R(Cone(f))$. Then the class $[Cone(f),triv]\in K_1(R,\Sigma_T)$ corresponds to $Det_{R_T}(f_T)^{-1}\in K_1(R_T)$ (both are indeed described by the same
endomorphism of $P_T^+\oplus P_T^-$ coming from $can_{Cone(f)_T}$).

\smallskip



Crucial to us will be the following:

\begin{lem} \label{frac} Consider morphisms $a,b: C\rightarrow C'$ in $D^p(R)$ whose localizations $a_T,b_T$ are isomorphisms in $D^p(R_T)$.
In $K_0(\Sigma_T)$, one has the equality $$\partial Det_{R_T}(a_Tb_T^{-1})=[Cone(a)]-[Cone(b)]$$ as long as the following condition holds: \smallskip

\noindent  $(frac)$ In $D^p(R)$, there exists a commutative square of the form $$\diagram{C&\hfl{a}{}&C'\cr \vfl{b}{}&&\vfl{d}{}\cr
C'&\hfl{c}{}&C'}$$ in which $c$ and $d$ also become isomorphisms after localization by $T$.
\end{lem}
Proof. An easy diagram chasing shows that $[Cone(c)]-[Cone(a)]=[Cone(d)]-[Cone(b)]$. Since $Det_{R_T}(a_Tb_T^{-1})=Det_{R_T}(d_T)^{-1}Det_{R_T}(c_T)$ and
$\partial$ is a homomorphism, it thus suffices to prove that $\partial Det_{R_T}(c_T)=[Cone(c)]$ and $\partial Det_{R_T}(d_T)=[Cone(d)]$. But this is clear
from the alternative characterization of the isomorphism $K_1(R,\Sigma_T)\simeq K_1(R_T)$ mentioned above.
\begin{flushright}$\square$\end{flushright}

We do not know whether or not the condition $(frac)$ always holds, neither if this computation of $\partial Det_{R_T}(a_Tb_T^{-1})$ is always correct.
We thus content ourselves with the following sufficient condition.

\begin{lem} \label{exfrac} Assume $R$ is Noetherian and consider $a,b: C\rightarrow C'$ in $D^p(R)$. If there exists $f$
in the center of $R$ such that $R[{1\over f}]\otimes^L_R Cone(b)=0$, then one can find $c,d:C'\rightarrow C'$ such that $da=cb$. In fact one can chose
$d=f^n$ for $n$ large enough.
\end{lem}
Proof. Since cohomology modules of $Cone(b)$ are finitely generated and almost all of them are zero, it is possible to find $n$ such that
$f^n:Cone(b)\rightarrow Cone(b)$ is zero. But then, the long exact sequence of $Ext$'s $$\cdot\cdot\cdot\, \rightarrow
Hom(C',C')\mathop\rightarrow\limits^{Hom(b,C')} Hom(C,C')\rightarrow Ext^1(Cone(b),C')\rightarrow\, \cdot\cdot\cdot$$ shows that $f^na$ has a trivial
image in $Ext^1(Cone(b),C')$ and thus comes from some $c\in Hom(C',C')$ as claimed.
\begin{flushright}$\square$\end{flushright}

\subsection{The evaluation map at Artin representations} \label{defevrho} ~~ \smallskip \\


Consider a profinite group $G=\limp_n G_n$ and a closed normal subgroup $H$ such that $\G:=G/H$ is isomorphic to $\Zp$. We use the following notations.

- We let $\La=\La(G):=\ilim \Z_p[G_n]$ denotes the Iwasawa algebra of $G$. We assume that this ring is left  Noetherian and regular (which will be the case when considering the Galois group of a $p$-adic Lie extension as in the introduction).

- If $O$ is the ring of integers of a finite extension $L$ of $\Q_p$, then $\La_O(G):=\ilim O[G_n]\simeq O\otimes_\Zp\La(G)$ and
$\La_L(G):=L\otimes_O\La_O(G)$ have similar properties.

- $S$ and $S^*$ denote the canonical Ore sets defined in \cite{CFKSV}. Recall that an element $f$ of $\La(G)$ is in $S$ if and only if $\La(G)/f$ is
a finitely generated $\La(H)$-module whereas $S^*=\cup_k p^kS$. As usual, $\Lambda_S=\Lambda(G)_{S}$ and $\Lambda_{S^*}=\Lambda(G)_{S^*}$ denote the corresponding localizations of $\L$. If $G=\G$ then $Q_O(\G):=\Lambda(\G)_{S^*}$ is the fraction field of $\Lambda_O(\G)$.

-$\mathfrak{M}_H(G)$ denotes the category of $S^*$-torsion finitely generated $\La(G)$-modules. We recall that a finitely generated $\La(G)$ module
$M$ is $S$-torsion (resp. $S^*$-torsion) if and only if it is a finitely generated $\La(H)$-module (resp. modulo its $p$-torsion).
In this context, the localization exact sequence reads:

$$\diagram{K_1(\La)&\hfl{}{}& K_1(\La_{S^*})&\hfl{\partial}{}& K_0(\mathfrak M_H(G))\hfl{}{}& K_0(\La)}$$

- If $R$ is one of the previous (localized) Iwasawa algebras, one often prefers to endow the dual of a left module with a left action, deduced from the right
one via the involution $g\mapsto g^{-1}$. This is our convention for duality of normic systems and their limit modules.
In this paragraph though, we leave right modules on the right, for the sake of clarity.

Consider a free $O$-module $V$ of finite rank. For any $O$-algebra $R$
(eg. $R=\La_O(G)$), we denote $V_R$ (resp. ${_RV}$) the right (resp. left) $R$-module $V\otimes_O R$ (resp. $R\otimes_O V$)
deduced from $V$ by extending scalars from $O$ to $R$. Also, $({_RV})^*:=Hom_R({_RV},R)$  (resp. $({V_R})^*:=Hom_R({V_R},R)$) is systematically given
its right (resp. left) $R$-module  structure coming from the $(R,R)$-bimodule structure of $R$.

It might be useful to recall that one has a canonical isomorphism of left (resp. right) $R$-modules   $(V_R)^*\simeq {_R(V^*)}$ (resp. $({_RV})^*\simeq (V^*)_R$) and canonical isomorphisms of $O$-algebras
$End_R((V_R)^*)\simeq End_R(V_R)^{op}$ and $End(({_RV})^*)\simeq End({_RV})^{op}$ ($R$-algebras). If now $R$ is a central $O$-algebra (eg. if $R=\L_O(G)$) we have moreover canonical isomorphisms of $O$-algebras $End_R({_RV})\simeq
R^{op}\otimes_O End_O(V)$, $End_R(V_R)\simeq End_O(V)\otimes_O R$. These isomorphisms are
subject to natural compatibilities, such as the commutativity of the following diagram:
$$\diagram{End_R(({_RV})^*)&\simeq& End_R({_RV})^{op}&\simeq &(R^{op}\otimes_O End_O(V))^{op}\cr |\wr&&&&|\wr\cr End_R((V^*)_R)&\simeq &
End_O(V^*)\otimes_OR&\simeq &End_O(V)^{op}\otimes_OR}$$

\medskip

Consider now an $O$-Artin representation $\rho:G\rightarrow Aut_O(V)$. By $O$-Artin representation we mean that $V$ is as above and $\rho$ has a finite image. Consider the unique homomorphism of $\Zp$-algebras $$\Phi:\La(G)\rightarrow End_{\La_O(G)}(V_{\La_O(G)})$$
sending $g\in G\subset \La_O(G)$ to $v\otimes \lambda\mapsto \rho(g)(v)\otimes g\lambda$.


By functoriality of $K_1$, one has a commutative diagram:
$$\xymatrix{K_1(\La(G))\ar @{->} @<+0pt> `u[r] `[rrr]^-{tw_\rho} [rrr]
\ar[rrd]_{K_1(\rho)}\ar[rr]^{K_1(\Phi)\hspace{.8cm}}&\hspace{.5cm}&K_1(End_{\La_O(G)}(V_{\La_O(G)}))\ar[d]^{K_1(\epsilon_V)}
\ar[r]_{\hspace{1cm}\sim}^{\hspace{1cm}Morita}& K_1(\La_O(G))\ar[d]^{K_1(\epsilon)}\cr
&&K_1(End_O(V))\ar[r]_{\hspace{.8cm}\sim}^{\hspace{.8cm}Morita}& K_1(O)}$$ where $tw_\rho$ is the map sending $[P,f]\in K_1(\La(G))$ to
$[(V_{\La_O(G)})^*\otimes_{\La(G)} P,id\otimes f]\in K_1(\La_O(G))$ where $(V_{\La_O(G)})^*$ is viewed as a $(\La_O(G),\La(G))$-bimodule, the left action of $\La_O(G)$ being the obvious one and the right action of  $\Lambda(G)$ being deduced via $\Phi$ from the natural right action of $End_{\La_O(G)}(V_{\La_O(G)})$.
Furthermore the image of this element in $K_1(O)$ is $[V^*\otimes_{\La(G)}P,id_{V^*}\otimes f]$ where the right $\La(G)$-module structure of $V^*$ is deduced from the right action of $G$ on $V^*$: $g\mapsto \rho(g)^*$. The composed homomorphism
$$\begin{array}{ccc}K_1(\La(G))&\rightarrow& K_1(O)\simeq O
^\times\cr
[P,f]&\mapsto& det_O(id_{V^*}\otimes f|V^*\otimes_{\La(G)}P)\end{array}$$ where $det_O$ denotes the usual (commutative) determinant, will simply be
denoted $\rho$.\smallskip

\begin{rem} \label{remtwist}
%
%
%
$(i)$ The element $tw_\rho([P,f])$ is also equal to $[V^*\otimes_\Zp
P,id\otimes f]$. Indeed, there are isomorphisms of left $\L_O(G)$ modules $$\begin{array}{rcl} (V_{\La_O(G)})^*\otimes_{\La(G)} P&\mathop\leftarrow\limits^\sim &_{\La_O(G)}(V^*)\otimes_{\La(G)} P\\ &\mathop\rightarrow\limits^\sim & V^*\otimes_\Zp P\end{array}$$
where the $(\L_O(G),\L(G))$-bimodule structure of $_{\La_O(G)}(V^*)$ is given by $h.(\lambda\otimes \phi)=h\lambda\otimes \phi$ and $(\lambda\otimes \phi).h=\lambda h\otimes \phi\circ (\rho(h))$), the left action of $\L_O(G)$ on $V^*\otimes_\Zp P$ is given by $h.(\phi\otimes p)=\phi\circ (\rho(h^{-1}))\otimes hp$, the first isomorphism is induced by $_{\La_O(G)}(V^*)\to (V_{\L_O(G)})^*$, $\lambda\otimes \phi\mapsto (v\otimes \mu\mapsto \lambda \phi(v)\mu)$ and the second is given by $g\otimes \phi\otimes p\mapsto \phi\circ \rho(g^{-1})\otimes gp$.

$(ii)$ Let us examine the simple case where $P$ is free with basis $e_1,\dots, e_n$ to fix the ideas. The element $[V^*\otimes_\Zp
P,id\otimes f]$ occurring in $(i)$ above has then a more convenient description as follows. Sending $\phi\otimes ge_i$ to $\phi\circ \rho(g)\otimes ge_i$ realises an isomorphism of $V^*\otimes_\Zp P$ with $V^*\otimes_\Zp P$ viewed as a left $\L_O(G)$-module via $h.(\phi\otimes \lambda)=\phi\otimes h\lambda$. Through this isomorphism $id\otimes f$ translates as the automorphism sending $\phi\otimes \lambda_je_j$ to $\iota(\lambda_j f(e_j))$ where $\iota:V^*\otimes_\Zp P\to V^*\otimes_\Zp P$ sends $\phi\otimes ge_i$ to $\phi\circ \rho(g)\otimes ge_i$. In particular we find that via the determinant isomorphism $K_1(O)\simeq O^\times$, $\rho([P,f])$ is nothing but the determinant of the automorphism of $O\otimes_{\L_O(G)}(V^*\otimes_\Zp P)\simeq \oplus_i V^*\otimes e_i$ sending $\phi\otimes e_j$ to $\sum_i\phi\circ \rho(f_{i}(e_j))$ if $f(e_j)=\sum f_i(e_j)e_i$.
\end{rem}

Next, following \cite{CFKSV}, we consider the extension of $\rho$ to a map \begin{eqnarray} \label{defrho} \rho:K_1(\La(G)_{S^*})\rightarrow L^\times \cup\{0;\infty\}\end{eqnarray}  defined as $\epsilon \rho_{Q_O(\G)}$ where:

- the map $\epsilon:Q_O(\G)^\times \rightarrow L\cup\{\infty\}$ coincides with $\L_O(\G)_I\to L$, the localization of the augmentation map at the augmentation ideal $I$, and takes the value $\infty$ elsewhere (ie. at the elements of $Q_O(\G)$ which are not integral at $I$).

- the map $\rho_{Q_O(\G)}:K_1(\La(G)_{S^*})\to Q_O(\G)^\times$ is defined by composing the obvious localized version of $tw_\rho$, $K_1(\L(G)_{S^*})\to K_1(\L_O(G)_{S^*})$, the fonctioriality map $K_1(\L_O(G)_{S^*})\to K_1(Q_O(\G))$ and the isomorphism $det_{Q_O(\G)}:K_1(Q_O(\G))\simeq Q_O(\G)^\times$.
\medskip

The map $\rho$ is multiplicative in the sense of the usual partial multiplication  $$(L\times L)\cup
((L^\times\cup\{\infty\})\times (L^\times\cup\{\infty\}))\rightarrow L\cup\{\infty\}$$
This means that for $x,y\in K_1(\La(G)_{S^*})$:

$(i)$  $\rho(x^{-1})=\rho(x)^{-1}$.

$(ii)$ the formula $\rho(x)\rho(y)=\rho(xy)$ is true as soon as it makes sense, ie. whenever $(\rho(x),\rho(y))\ne (0,\infty), (\infty,0)$.
\medskip

\begin{lem} \label{rhodual} Let $[P,f]$ in $K_1(\Lambda(G)_{S^*})$ and consider its dual $[P^\vee,f^\vee]$ where $(-)^\vee$ stands for $Hom_{\Lambda(G)_{S^*}}(-,\Lambda(G)_{S^*})$  viewed as a left $\L(G)_{S^*}$-module via the involution $g\mapsto g^{-1}$. If $\rho$ is an Artin representation with contragredient $\rho^\vee$ then in $K_1(L)$: $$\rho([P,f])=\rho^\vee([P^\vee,f^\vee])$$
\end{lem}

Proof. One easily reduces to the analogous statement for $\L=\L(G)$ in place of $\Lambda(G)_{S^*}$. In virtue of Rem. \ref{remtwist} $(i)$, it is then sufficient to prove that for a finitely generated projective $P$ over $\L(G)$ we have in $L^\times$:  $$det_L(id\otimes id\otimes f|L\otimes_{\L_O}(V^*\otimes_{\Zp} P))=det_L(id\otimes id\otimes Hom_{\L}(f,\L)|L\otimes_{\L_O}((V^\vee)^*\otimes_\Zp Hom_\L(P,\L)))$$
the action of $\L_O$ on $V^*\otimes_\Zp P$ (resp. $(V^\vee)^*\otimes_\Zp Hom_\L(P,\L)=V\otimes_\Zp Hom_\L(P,\L)$) being given by $g.(\phi\otimes p)=\phi\circ \rho(g^{-1})\otimes gp$ (resp. $g.(v\otimes \phi)=\rho(g)v\otimes (p\mapsto \phi(p)g^{-1})$). Since the determinant of an endomorphism and its transpose are equal it is sufficient to build an isomorphism \begin{eqnarray} \label{isovoulu} Hom_L(L\otimes_{\L_O}(V^*\otimes_{\Zp} P),L)\simeq L\otimes_{\L_O}((V^\vee)^*\otimes_\Zp Hom_\L(P,\L))\end{eqnarray} identifying $Hom_L(id\otimes id\otimes f,L)$ to $id\otimes id\otimes Hom_{\L}(f,\L)$. Let us examine the left term of (\ref{isovoulu}). Let us chose $n$ such that $\rho$ factors through $G_n=G/H_n$. Denoting $V_L=L\otimes_O V$ we have the following series of natural isomorphisms $$\begin{array}{rcl}Hom_L(L\otimes_{\L_O}(V^*\otimes_{\Zp} P),L)&\simeq & Hom_L((V_L^*\otimes_{\Zp} P)_G,L)\\ &\simeq & Hom_L(V_L^*\otimes_\Zp P_{H_n},L)^{G_n} \\ &\simeq & (V_L\otimes_L Hom_\Zp(P_{H_n},L))^{G_n} \\
&\simeq &  (V_L\otimes_L Hom_\Zp(P_{H_n},L))_{G_n}\end{array}$$
where the invariants or coinvariants occurring respectively in the second, third, third and fourth term  are taken with respect to the following left actions of $G$ or $G_n$: $g.(\phi\otimes p)=\phi\circ \rho(g^{-1})\otimes gp$, $g.\psi=(\phi\otimes p\mapsto \psi(\phi \circ \rho(g)\otimes g^{-1}p))$, $g.(v\otimes \phi)=\rho(g)v\otimes (p\mapsto \phi(g^{-1}p))$,  $g.(v\otimes \phi)=\rho(g)v\otimes (p\mapsto \phi(g^{-1}p))$. Regarding the second term of (\ref{isovoulu}) we have $$\begin{array}{rcl} L\otimes_{\L_O}((V^\vee)^*\otimes_\Zp Hom_\L(P,\L))&\simeq & (V_L\otimes_\Zp Hom_\L(P,\L))_G\\ &\simeq & (V_L\otimes_\Zp Hom_\L(P,\L)_{H_n})_{G_n} \\ &\simeq & (V_L\otimes_L Hom_{\Zp[G_n]}(P_{H_n},L[G_n]))_{G_n}\end{array}$$
where the coinvariants occurring respectively in the second, third and fourth term are taken with respect to the following left action of $G$, $H_n$, $G_n$ and $G_n$: $g.(v\otimes \phi)=\rho(g)v\otimes (p\mapsto \phi(p)g^{-1})$, $h.\phi=(p\mapsto \phi(p)h^{-1})$, $g.(v\otimes \phi)=\rho(g)v\otimes (p\mapsto \phi(p)g^{-1})$, $g.(v\otimes \phi)=\rho(g)v\otimes (p\mapsto \phi(p)g^{-1})$. We finally obtain (\ref{isovoulu}) by noticing the isomorphism $$\begin{array}{rcl}  Hom_\Zp(P_{H_n},L) &\simeq & Hom_{\Zp[G_n]}(P_{H_n},L[G_n]) \\ \phi &\mapsto & (p\mapsto \sum_{g\in G_n} \phi(gp)[g^{-1}])\end{array}$$
The reader may easily check by himself that (\ref{isovoulu}) is compatible with the automorphisms induced by $f$ as desired.
\begin{flushright}$\square$\end{flushright}

The next results concerning the computation of (\ref{defrho}) composed with (\ref{defDet}) will be useful to us.

\begin{lem}\label{rho_Det} Consider an automorphism $f$ of an object $C$ of $D^p(\La(G)_{S^*})$ as well as its determinant $Det_{\La(G)_{S^*}}(f)\in
K_1(\La(G)_{S^*})$ (cf. (\ref{defDet})). In $Q_O(\G)^\times$ we have:
$$\rho_{Q_O(\G)}(Det_{\La(G)_{S^*}}(f))=
\prod_qdet_{Q_O(\G)}(H^q(id\otimes f)|H^q((V_{Q_O(\G)})^*\Ltens_{\La(G)_S^*}C))^{(-1)^q}$$
where the tensor product occurring on the right hand side is taken with respect to the right $\La(G)_{S^*}$-module structure of $(V_{Q_O(\G)})^*$ defined using $\Phi$ in the obvious way.
\end{lem}
Proof. We can always assume that $C$ is strictly perfect. In this case we have by definition of $\rho_{Q_O(\G)}$ and compatibility of $Det$ with scalar extension:
$$\begin{array}{rcl}\rho_{Q_O(\G)}(Det_{\La(G)_{S^*}}(f))&=&det_{Q_O(\G)}(Det_{Q_O(\G)}(id\otimes f))\end{array}$$ where $id\otimes f$ is
the automorphism of $(V_{Q_O(\G)})^*\otimes_{\La(G)_S^*}C\in D^p(Q_O(\G))$ deduced from $f$. Now $Q_O(\G)$ being a field, the
complex $(V_{Q_O(\G)})^*\otimes_{\La(G)_{S^*}}C$ is cohomologically perfect, and thus
$$\begin{array}{rcl}det_{Q_O(\G)}Det_{Q_O(\G)}(id\otimes f)&=&
\prod_qdet_{Q_O(\G)}Det_{Q_O(\G)}(H^q(id\otimes f))^{(-1)^q}\\
&=& \prod_qdet_{Q_O(\G)}(H^q(id\otimes f)|H^q((V_{Q_O(\G)})^*\otimes_{\La(G)_{S^*}}C))^{(-1)^q}\end{array}$$
\begin{flushright}$\square$\end{flushright}

\begin{lem}\label{Tor} Let $\epsilon$, $I$, $\La_O(\G)_I$ as in (\ref{defrho}).
Consider a finitely generated $\La_O(\G)_I$-module $M$ together with an endomorphism $f:M\rightarrow M$ such that $Q_O(\G)\otimes_{\La_O(\G)_I}
f$ is invertible. The formula
$$\begin{array}{l}\epsilon det_{Q_O(\G)}(id\otimes f|Q_O(\G)\otimes_{\La_O(\G)_I}M)\hspace{5cm}\\ \hspace{1cm}=
det_L(Tor_0^{\La_O(\G)_I}(id,f)|Tor_0^{\La_O(\G)_I}(L,M)).det_L(Tor_1^{\La_O(\G)_I}(id,f)|Tor_1^{\La_O(\G)_I}(L,M))^{-1}\end{array}$$ is true
whenever the right hand term makes sense in $L^\times\cup\{0,\infty\}$. \end{lem} Proof. Let $M'\subset M$ denote the $\La_O(\G)_I$-torsion submodule of $M$,
and $M'':=M/M'$. Let us denote $f':M'\rightarrow M'$, $f'':M''\rightarrow M''$ the endomorphisms induced
by $f$.

On the one hand, $Q_O(\G)\otimes_{\La_O(\G)_I}M'=0$, and thus $$det_{Q_O(\G)}(id\otimes f|Q_O(\G)\otimes_{\La_O(\G)_I}M)=det_{Q_O(\G)}(id\otimes
f''|Q_O(\G)\otimes_{\La_O(\G)_I}M'')$$

On the other hand, $M''$ is free since $\La_O(\G)_I$ is a discrete valuation ring. We thus have  a short exact sequence and an isomorphism
$$\diagram{Tor_0^{\La_O(\G)_I}(L,M')&\hookrightarrow &Tor_0^{\La_O(\G)_I}(L,M)&\twoheadrightarrow &Tor_0^{\La_O(\G)_I}(L,M'')\cr
Tor_1^{\La_O(\G)_I}(L,M')&\simeq &Tor_1^{\La_O(\G)_I}(L,M)&&}$$ both of which are compatible with $f$. Whence equalities:
$$\begin{array}{l}det_L(Tor_0^{\La_O(\G)_I}(id,f)|Tor_0^{\La_O(\G)_I}(L,M))\hspace{5cm}\\ \hspace{1cm}=
det_L(Tor_0^{\La_O(\G)_I}(id,f)|Tor_0^{\La_O(\G)_I}(L,M'))det_L(Tor_0^{\La_O(\G)_I}(id,f)|Tor_0^{\La_O(\G)_I}(L,M''))
\end{array}$$
$$det_L(Tor_1^{\La_O(\G)_I}(id,f)|Tor_1^{\La_O(\G)_I}(L,M))=det_L(Tor_1^{\La_O(\G)_I}(id,f)|Tor_1^{\La_O(\G)_I}(L,M'))$$
Now $M'$ is of finite length, and is thus subject to Koszul duality: $$Hom_L(Tor_1^{\La_O(\G)_I}(L,M'),L)\simeq Tor_0^{\La_O(\G)_I}(L,Hom_L(M',L))$$
Since $det_L$ isn't affected by $L$-linear duality, this shows that
$$det_L(Tor_1^{\La_O(\G)_I}(id,f)|Tor_1^{\La_O(\G)_I}(L,M'))=det_L(Tor_0^{\La_O(\G)_I}(id,f)|
Tor_0^{\La_O(\G)_I}(L,M'))$$

Now, in the formula to be proven, we see that the right hand term makes sense if and only if
$det_L(Tor_0^{\La_O(\G)_I}(id,f)|Tor_0^{\La_O(\G)_I}(L,M'))$ is non zero (ie. if and only if \break $det_L(Tor_1^{\La_O(\G)_I}(id,f)|Tor_0^{\La_O(\G)_I}(L,M))$ is non zero) in which case the factors \break $det_L(Tor_q^{\La_O(\G)_I}(id,f)|Tor_0^{\La_O(\G)_I}(L,M'))$, $q=0,1$, cancel each other and the desired formula reduces to the obvious equality $$\epsilon
det_{Q_O(\G)}(id\otimes f|Q_O(\G)\otimes_{\La_O(\G)_I}M'')=det_L(Tor_0^{\La_O(\G)_I}(id,f)|Tor_0^{\La_O(\G)_I}(L,M''))$$
\begin{flushright}$\square$\end{flushright}

%
%
%

\section{Selmer complexes and crystalline cohomology}
\grand

In this section we define Selmer complexes for abelian varieties over a one variable function field. We begin with their basic properties and a duality theorem. In the semistable case we pursue with a review of the main result of \cite{TV}, which will be the cornerstone of our proof for the Iwasawa main conjecture.

\subsection{Preliminaries}
Let us begin with some technical facts regarding derived categories.

\para \label{prelimDI} We will frequently use derived categories of the type $D(\cal C^{\cal I})$, where $\cal C^{\cal I}$ denotes the category of contravariant functors $\cal I\rightarrow \cal C$, $\cal C$ being an abelian category. Note that $\cal C^\I$ is an abelian category. If $\cal C$ has enough injectives and it has products indexed by the subsets of $Ob(\I)$, then the same holds for $\cal C^\I$.
If $\cal C$ is the category of modules
on a ringed topos $(E,A)$, $Mod(E,A)^\I$ naturally identifies with the category of modules on the topos $(E^\I,A)$. This point of view offers the possibility to consider the more general category $Mod(E^{\cal I},A_.)$, with $A_.$  a projective system of rings of $E$ indexed by $\cal I$.

These categories are especially useful for some particular choices of $\I$ which we explain now.

- Let $\cal I=\mathbb N$ viewed as the category where $Hom(k,k')$ has one element if $k\le k'$ and is empty otherwise (resp. $\I=\N^{op}$, resp. $\cal I=N_2:=\N\sqcup_{|\N|}\N^{op}$). In that case, the objects of the category $\cal C^\I$ are the projective systems of $\cal C$ indexed by integers (resp. the inductive systems of $\mathcal C$ indexed by integers, resp. the triples $((A_k)_k,(p_{k',k})_{k'\ge k},(i_{k,k'})_{k\le k'})$ where $((A_k)_k,(p_{k',k})_{k'\ge k})$ is in $\cal C^\N$ and  $((A_k)_k,(i_{k,k'})_{k\le k'})$ is in $\cal C^{\N^{op}}$).

- Let $\cal I=[1]^a$ for some integer $a\ge 1$, $[1]$ denoting the category $\{0,1,\le\}$. In that case, we think of objects of  $D(\cal C^{\cal I})$ as $(a+1)$-uple naive complexes which are zero outside a specified range of the form  $]-\infty,+\infty[\times [i_1,i_1+1]\times \dots, \times [i_a,i_a+1]$. The interest of this category lies in the fact that forming total complexes (with an appropriate sign convention which the interested reader is invited to specify) give rise to a triangulated functor $Tot_{a+1}:D(\cal C^{[1]^a})\rightarrow D(\cal C)$.

\para \label{prelimMF} The case $\cal I=[1]$ is already interesting since $Tot_2$ gives rise to a functorial version of the mapping cone / mapping fiber construction. Let us discuss a fact that will be used repeatedly. A natural transformation $t:F_1\rightarrow F_2$ between functors $\cal C\rightarrow \cal C'$ can be thought as a functor $F$ taking its values in $\cal C'^{[1]}$. If $\cal C$ has enough injectives we can form the derived functor $RF:D(\cal C)\rightarrow D(\cal C'^{[1]})$. Next, we define $MF(t):=Tot_2\circ RF$. Since the mapping fiber construction gives a functor from $D(\cal C'^{[1]})$  to that of distinguished triangles in $D(\cal C')$ we find a canonical distinguished triangle $MF(t)(C)\rightarrow RF_1(C)\rightarrow RF_2(C)\rightarrow MF(t)(C)[1]$ varying functorially with respect to $C$ in $D(\cal C)$. Consider furthermore $F_0(C):=Ker(F_1(C)\rightarrow F_2(C))$. Then $F_0$ defines a right derivable functor. Since $RF_0\simeq MF(t_0)$ for $t_0:F_0\rightarrow 0$ and $t_0$ canonically maps to $t$, one gets a canonical morphism $RF_0\rightarrow MF(t)$. This natural transformation is an isomorphism as long as the morphism $F_1(C)\rightarrow F_2(C)$ is epimorphic for injective $C$ (because in that case the sequence $0\rightarrow F_0(C)\rightarrow F_1(C)\rightarrow F_2(C)\rightarrow 0$ is exact and the natural transformation in question is the one occurring in the construction of the distinguished triangle associated to a short exact sequence).

 \petit

\grand
\subsection{Flat cohomology of $C$ vanishing at $Z$} ~~ \\

We review the definition of the functor of vanishing cohomology, which is a necessary tool to even state the comparison result we need from \cite{TV} (see Sect. \ref{reminderTV}).



\para If $S$ is any scheme, we denote $$\epsilon:S_{FL}\rightarrow S_{et}$$ the natural morphism from the big flat topos of $S$ (ie. the category of sheaves on $Sch/S$ endowed with the topology generated by surjective families of flat morphisms of finite type) to the small \' etale one (ie. the category of sheaves one the category $Et/S$ of \' etale $S$-schemes endowed with the topology generated by surjective families, cf SGA 4, VII, 1.2). The morphism $\epsilon$ is functorial with respect to $S$, ie. there is a canonical isomorphism $f\epsilon_S \simeq \epsilon_{S'} f$ for any $f:S_{et}\rightarrow S'_{et}$ induced by a morphism of schemes $f:S\rightarrow S'$. However, the reader should be aware that the base change morphism $f^{-1}\epsilon_*\rightarrow \epsilon_*f^{-1}$ is certainly not an isomorphism in general if $f$ is not \' etale.

\begin{defn} \label{defvanishingcoh} Keep the above notations and let furthermore $i:T\hookrightarrow S$ be the inclusion of a closed subscheme as well as $j:Y\rightarrow S$ that of its open complement. We also assume given a projective system of torsors $\underline S=(S_n,G_n)$ where $G_n$ is a finite quotient of a profinite group $G$ and denote $\G(\underline S,-):Mod(S_{et},\Z)\rightarrow {_{\underline G}Mod(\Z)}$ the associated \emph{normic sections functor}.  \petit

We define the following.  \debrom
\item The functor of \emph{global sections vanishing at $T$}, $$\G^T(S,-):Mod(S_{FL},\Z)\rightarrow Mod(\Z)$$ is defined as $\G^T(S,F):=Ker(F(S)\rightarrow F(T))$.
\item The functor of \emph{\' etale sections vanishing at $T$}, $$\epsilon^T_*:Mod(S_{FL},\Z)\rightarrow Mod(S_{et},\Z)$$ is defined as $\epsilon^T_*F:=Ker(\epsilon_*F\rightarrow \epsilon_*i_*i^{-1}F)$.
\item The functor of \emph{normic sections vanishing at $T$}, $$\G^T(\underline S,-):Mod(S_{FL},\Z)\rightarrow {_{\underline G}Mod(\Z)}$$ is defined as $\G^T(\underline S,F):=\G(\underline S,\epsilon^T_*F)$. Its $n^{th}$ component is thus $Ker(F(S_n)\rightarrow F(T_n))$ where $T_n:=S_n\times_ST$.
\finrom

\end{defn}

Note that by definition, one may respectively retrieve $\G^T(S,-)$ and $\G^T(\underline S,-)$ from $\epsilon^T_*$ by forming global and normic sections (Def. \ref{defGnorm}).

\begin{lem} \label{tdT} The above functors are right derivable and their derived functors are subject to natural  distinguished triangles in $D(S_{et})$, $D(Mod(\Z))$ and $D({_{\underline G}Mod(\Z)})$ which are functorial with respect to $F$ in $D^+(S_{FL})$:
$$\diagram{
R\epsilon^T_*F&\hfl{}{}&R\epsilon_*F&\hfl{}{}&Ri_*R\epsilon_*F_{|T}&\hfl{+1}{}\cr
R\G^T(S,F)&\hfl{}{}&R\G(S,F)&\hfl{}{}&R\G(T,F_{|T})&\hfl{+1}{}\cr
R\G^T(\underline S,F)&\hfl{}{}&R\G(\underline S,F)&\hfl{}{}&R\G(\underline T,F_{|T})&\hfl{+1}{}
}$$
where $\underline T$ and $\underline Y$ denote the projective systems of torsors deduced from $\underline S$. Taking $n$th components in the last distinguished triangles gives
$$\diagram{ R\G^{T_n}(S,F)&\hfl{}{}&R\G(S_n,F_{|S_n})&\hfl{}{}&R\G(T_n,F_{|T_n})&\hfl{+1}{}}$$
\end{lem}
Proof. Since injective abelian sheaves of $S_{FL}$ are flasque in the sense of [SGA 4, exp V, 4.7], it follows from the discussion \ref{prelimMF} that $R\epsilon^T_*$ identifies with the functorial mapping fiber of the natural transformation $R\epsilon_*\rightarrow R(\epsilon_*i_*i^{-1})$ arising from $\epsilon_*\rightarrow \epsilon_*i_*i^{-1}$. Now, $i^{-1}$ preserves injectives ($i$ is a localization morphism) and thus $ R(\epsilon_*i_*i^{-1})\simeq R(\epsilon_*i_*)i^{-1}$. The first distinguished triangle follows since $R(\epsilon_*i_*)\simeq Ri_*R\epsilon_*$. The other cases are similar.\begin{flushright}$\square$\end{flushright}

\petit

\grand

\subsection{Selmer complexes} ~~ \\

The purpose of this section is to give a tractable definition for (normic) Selmer complexes, to establish their relation to the complex of (normic) cohomology vanishing at $Z$ (the relation with the complexes $R\G_{ar,V}$ appearing in \cite{KT} will also be explained briefly in Rem. \ref{SelGar}), as well as a duality theorem. \petit

\para \label{intro2}
We consider the following situation:

\begin{eqnarray}\label{diagcomplet}
\diagram{Z_v&\hfl{z_v}{}&C_v&\hflrev{j_v}{}&U_v\cr  \vfl{{\iota_{Z_v}}}{}&&\vfl{{\iota_{C_v}}}{}&&\vfl{{\iota_{U_v}}}{}\cr
Z&\hfl{z}{}&C&\hflrev{j}{}&U\cr}
\end{eqnarray}

where: \\
- $C$ is a connected proper smooth $1$-dimensional $\Fp$-scheme with function field $K$ and constant field $k$. \\
- $Z$ is an effective divisor on $C$, ie. a finite union of $0$-dimensionnal irreducible closed subschemes $Z_v$. We denote $Z_v^{red}=Spec(k_v)$ and $Z^{red}=\sqcup_v Z_v^{red}$ the underlying reduced schemes. \\
- $U$ denotes the open subcheme of $C$ which is complementary to $Z$. \\
- For each $v$ in $|Z|$, we denote $O_v$ the completion of the local ring of $C$ at $v$, $K_v$ its fraction field,  $C_v=Spec(O_v)$ and $U_v=Spec(K_v)$. \\
- The arrows $i,j,j_v,i_v,{\iota_{U_v}},{\iota_{C_v}}$ and $
{\iota_{Z_v}}$ are the obvious ones. \petit

We fix moreover a profinite Galois extension  $K_\infty/K$ with group $G=\ilim  G_n$ satisfying the following properties: \\
- $G$ is a $p$-adic Lie group without $p$-torsion, \\
- $K_\infty/K$ is unramified everywhere and contains the constant $\Zp$-extension $K_{ar}=Kk_\infty/K$ whose Galois group is denoted $\Gamma$.

Thanks to the second assumption, there is an essentially unique projective system of torsors $\underline C=(C_n,G_n)$ over $C$ corresponding to the Galois tower $\underline K=(K_n,G_n)$ over $F$. \petit

Finally we fix an abelian variety $A$ over $K$ and denote $\cal A$ its N\'eron model over $C$. We always assume that $A$ has good reduction outside $Z$ (ie. $A_{|U}$ is an abelian scheme over $U$).

%
%
%
%

\para \label{auxSel}
We place ourselves in the situation described in \ref{intro2}. The divisor $Z$ will be allowed to vary, but its support will always contain the points of bad reduction of $A$. Unless specified otherwise, cohomology is meant in the sense of big flat topoi. 
To define the Selmer complex we need the following intermediary functors. We use the notations $\iota_{C_v}:C_v\rightarrow C$, $j_v:U_v\rightarrow C_v$ for all places of $C$ (not only those of $Z$).


- $Loc_{(U)}:Mod(C_{FL},\Z)\rightarrow Mod(C_{et},\Z)^{[1]}$ is the functor taking $F\in Mod(C_{FL},\Z)$  to
$$\left[\diagram{j_*\epsilon_*F_{|U}&\hfl{}{}&(\mathop\prod\limits_{v\in |U|} \iota_{C_v,*}\epsilon_*F_{|C_v})\times (\mathop\prod\limits_{v\in |Z|} \iota_{C_v,*}j_{v,*}\epsilon_*F_{|U_v})}\right]
$$

Let $\cal U$ denote the filtered set of open subschemes of $C$. Since for $Z\subset Z'$, there is an obvious natural transformation $Loc_{(U)}\rightarrow Loc_{(U')}$, we also have a functor $Loc_{(\cal U)}:Mod(C_{FL},\Z)\rightarrow Mod(C_{et},\Z)^{[1]\times \cal U}$.

- $\tau_{\ge[{0\, 1}]}:D(Mod(\Z)^{[1]})\rightarrow D(Mod(\Z)^{[1]})$ (resp. $D({_{\underline G}Mod(\Z)^{[1]}})\rightarrow D({_{\underline G}Mod(\Z)^{[1]}})$, resp.  $D({_{\underline G}Mod(\Z)^{[1]\times \cal U}})\rightarrow D({_{\underline G}Mod(\Z)^{[1]\times \cal U}})$) is described as $$\begin{array}{ccc}\left[\diagram{A&\hfl{}{}&B}\right]& \mapsto \left[\diagram{A&\hfl{}{}&\tau_{\ge 1}B}\right]\end{array}$$


\begin{defn} We define Selmer functors as follows.


$(i)$ The \emph{$(U)$-Selmer functor} $$Sel_{(U)}(C,-):D^+(C_{FL})\rightarrow D^+(Mod(\Z))$$ is defined as $Sel_{(U)}(C,-):=Tot_2\circ \tau_{\ge[{0\, 1}]}\circ R\G(C,-)\circ  RLoc_{(U)}$.

$(ii)$ The \emph{$(U)$-normic Selmer functor} $$Sel_{(U)}(\underline C,-):D^+(C_{FL})\rightarrow D^+({_{\underline G}Mod(\Z)})$$ is defined as $Sel_{(U)}(\underline C,-):=Tot_2\circ \tau_{\ge[{0\, 1}]}\circ R\G(\underline C,-)\circ  RLoc_{(U)}$, where $\G(\underline C,-)$ is the normic section functor of Def. \ref{defGnorm}.

$(iii)$ The \emph{$\cal U$-normic Selmer functor} $$Sel_{(\cal U)}(\underline C,-):D^+(C_{FL})\rightarrow D^+({_{\underline G}Mod(\Z)^{\cal U}})$$ is defined as $Sel_{(\cal U)}(\underline C,-):=Tot_2\circ \tau_{\ge[{0\, 1}]}\circ R\G(\underline C,-)\circ  RLoc_{(\cal U)}$.
\end{defn}

Let us point out that those Selmer functors are not triangulated, since their definition involves truncations.

\begin{rem} \label{remselmer} $(i)$ The functor $(ii)$ (resp. $(i)$) can be retrieved from functor $(iii)$ (resp. $(ii)$) by taking $U$-components of the inductive system (resp. component $0$ of the normic systems).

$(ii)$ It follows immediately from the definition that $Sel_{(U)}(\underline C,F)$ fits into a distinguished triangle  $$\diagram{Sel_{(U)}(\underline C,F)&\hflcourte{}{}&R\G(\underline U,F)&\hflcourte{}{}& (\mathop\prod\limits_{v\in |U|}\tau_{\ge 1}R\G(\underline C_v,F))\times  (\mathop\prod\limits_{v\in |Z|}\tau_{\ge 1}R\G(\underline U_v,F))&\hflcourte{+1}{}}$$
and similarly for $Sel_{(U)}(C,F)$. In the literature, Selmer complexes are usually designed to fit into a distinguished triangle where $\tau_{\ge 1}R\G(\underline C_v,F))$ does not appear.  This factor will disappear when $F$ is the N\'eron model of an abelian variety with good reduction over $U$ (see Prop. \ref{compsel} below).  The reason why  we add this term is to get functoriality with respect to $U$. Note that the functoriality of $Sel_{(U)}(\underline C,F)$ with respect to $U$ not only occurs in the derived category, but already at the level of complexes. This fact is reflected by the existence of $Sel_{(\cal U)}(\underline C,F)$.
\end{rem}





\begin{lem} \label{selindepU} Let $Z\subset Z'$, $U'=C-Z'$, $U=C-Z$, so that $U'\subset U$ and let $\cal A/C$ be as in \ref{intro2}.

$(i)$ The following natural morphisms in $D({_{\underline G}Mod(\Z)})$ are invertible:
$$\diagram{Sel_{(U)}(\underline C,\cal A)&\hfl{(1)}{}& Sel_{(U')}(\underline C,\cal A)&\hfl{(2)}{}& \limi_{\cal U} Sel_{(\cal U)}(\underline C,\cal A)}$$

$(ii)$ If $\Phi_v/Z_v^{red}$ denotes the component group of $\cal A$ at $v$, we have \begin{eqnarray}\label{sp}\begin{array}{rcl}
\tau_{\ge 1}  R\G(\underline C_v,{{\cal A}}_{|C_v})&\simeq & \tau_{\ge 1}R\G(\underline Z_v^{red},\cal A_{|Z_v^{red}})\\
&\simeq & H^1(\underline Z_v^{red},\Phi_v)[-1]\end{array}\end{eqnarray}
This group is in particular zero if $\cal A$ has good reduction at $v$.

\end{lem}
Proof. $(i)$ It suffices to prove the assertion about $(1)$. For any $F\in D^+(C_{FL})$, one has a tautological diagram in $D({_{\underline G}Mod(\Z)})$ with distinguished rows and columns as follows:
$$\diagram{Sel_{(U)}(\underline C,F)&\hfl{}{}&R\G(\underline U)&\hfl{}{}& (\mathop\prod\limits_{v\in |U|}\tau_{\ge 1}R\G(\underline C_v))\times  (\mathop\prod\limits_{v\in |Z|}\tau_{\ge 1}R\G(\underline U_v))&\hfl{+1}{}&\cr
\vflcourte{a}{}&&\vflcourte{b}{}&&\vflcourte{c}{}\cr
Sel_{(U')}(\underline C,F)&\hfl{}{}& R\G(\underline U')&\hfl{}{}& (\mathop\prod\limits_{v\in |U'|}\tau_{\ge 1}R\G(\underline C_v))\times  (\mathop\prod\limits_{v\in |Z'|}\tau_{\ge 1}R\G(\underline U_v))&\hfl{+1}{}&\cr
\vflcourte{}{}&&\vflcourte{}{}&&\vflcourte{}{}\cr
Cone(a)&\hfl{}{}&Cone(b)&\hfl{}{}&Cone(c)&\hfl{+1}{}& \cr
\vflcourte{+1}{}&&\vflcourte{+1}{}&&\vflcourte{+1}{}\cr}$$
where we wrote $R\G(\underline U)$ instead of $R\G(\underline U,F_{|U})$ by lack of space and similarly for $U'$, $C_v$, and $U_v$. The morphism $b$ is induced by restriction from $U$ to $U'$. The morphism $c$ is induced by restriction from $C_v$ to $U_v$ for $v\in |Z'|\backslash |Z|=|U|\backslash |U'|$ and is the identity for other $v$'s.  The definition of $a$, $Cone(a)$, $Cone(b)$, $Cone(c)$ is tautological. Let us compute these mapping cones.

To begin with, we note that by definition of relative cohomology
$$Cone(b)\simeq  R\G_{U\backslash U'}(\underline U, F_{|U})[1]$$
Next we note that there is an natural morphism $$\mathop\prod\limits_{v\in |Z'|\backslash |Z|} (R\G_{Z_v}(\underline{U}_v, F_{|U_v})[1])\to Cone(c)$$ which is invertible if the restriction morphisms $$R\G(\underline C_v,F_{|C_v})\rightarrow R\G(\underline U_v,F_{|U_v})$$ are invertible in degree $0$. This is the case for $F=\cal A$ by the N\'eron extension property.

The morphism $Cone(b)\rightarrow Cone(c)$ is thus induced by  $$R\G_{U\backslash U'}(\underline U, \cal A)\rightarrow  \mathop\prod\limits_{v\in |Z'|\backslash |Z|} R\G_{Z_v}(\underline{C}_v, \cal A)$$ and is thus an isomorphism by \cite{Mi1}, III, 1.28 and \cite{Mi2}, III, 7.14. It follows that $Cone(a)=0$, ie. that $a$ is an isomorphism.

$(ii)$ The first isomorphism is by \cite{Mi1}, III 3.11 and the second is by Lang's lemma for the smooth connected algebraic group (\cite{La}, Thm. 2) $\cal A_{|Z^{red}_{n,v}}^0$.
\begin{flushright}$\square$\end{flushright}

\para We are now in position to define the Selmer complex of $A/K$.
The reader will retrieve the usual description below (Prop. \ref{compsel} $(iii)$, Cor.  \ref{selcoh}).

\begin{defn} We define the following:

$(i)$ The \emph{normic Selmer complex of $A$} is defined in $D^b({_{\underline G}Mod(\Z)})$ as
$$Sel(A/\underline K):=\limi_{\cal U} Sel_{(\cal U)}(\underline C,\cal A).$$

$(ii)$ The \emph{normic $p^.$-Selmer complex of $A$} is defined in $D^b({_{\underline G}Mod(\Zp)^{N_2}})$ as $$Sel_{p^.}(A/\underline K):=\Z/p^.\Ltens Sel(A/\underline K)[-1]$$

$(iii)$ The normic \emph{$p^\infty$-Selmer complex of $A$} is defined in $D^b({_{\underline G}Mod(\Zp)})$ as    $$Sel_{p^\infty}(A/\underline K):=\limi Sel_{p^.}(A/\underline K)$$

$(iv)$ The normic \emph{$T_p$-Selmer complex of $A$} is defined in $D^b({_{\underline G}Mod(\Zp)})$ as    $$Sel_{T_p}(A/\underline K):=R\limp Sel_{p^.}(A/\underline K)$$
\end{defn}


The following proposition summarizes the relations between the normic Selmer complex of an abelian variety, its cohomology over $C$, its cohomology vanishing at $Z$, its cohomology over $U$ and its compactly supported cohomology over $U$.

\begin{prop} \label{compsel}

$(i)$ $Sel(A/\underline K)$ fits into canonical distinguished triangles as follows.
$$\diagram{Sel(A/\underline K)&\hfl{}{}R\G(\underline C,\cal A)&\hfl{}{}&\mathop\prod\limits_{v\in |Z|}H^1(\underline{C}_v,\cal A_{|C_v})[-1]&\hfl{+1}{}\cr
Sel(A/\underline K)&\hfl{}{}R\G(\underline U,\cal A_{|U})&\hfl{}{}&\mathop\prod\limits_{v\in |Z|}H^1(\underline{U}_v,\cal A_{|U_v})[-1]&\hfl{+1}{}\cr
Sel(A/\underline K)&\hfl{}{}R\G(\underline K,A)&\hfl{}{}&\mathop\oplus\limits_{v\in |C|}H^1(\underline{U}_v,A_{|U_v})[-1]&\hfl{+1}{}\cr}$$

$(ii)$ $Sel_{p^.}(A/\underline K)$ compares to cohomology over $C$ (resp. vanishing at $Z$) via canonical distinguished triangles as follows.
$$\begin{array}{c}Sel_{p^.}(A/\underline K)\rightarrow
R\G(\underline C,\Z/p^.\Ltens{\cal A}[-1])\rightarrow \mathop\prod\limits_{v\in |Z|}\Z/p^.\Ltens
H^1(\underline{k}_v,\Phi_v)[-2]\mathop\rightarrow\limits^{+1} \\
R\G^Z(\underline C,\Z/p^.\Ltens{\cal A}[-1])\rightarrow Sel_{p^.}(A/\underline K)\rightarrow \oplus_{v\in
|Z|}\Z/p^.\Ltens\cal A(\underline Z_v)[-1]\mathop\rightarrow\limits^{+1}\end{array}$$

$(iii)$ $Sel_{p^.}(A/\underline K)$ compares to (resp. compactly supported) cohomology over $U$ via canonical distinguished triangles as follows.
$$\diagram{Sel_{p^.}(A/\underline K)&\hfl{}{}&R\G(\underline U,({\cal A}_{|U,p^.}))&\hfl{}{}& \mathop\prod\limits_{v\in |Z|}\Z/p^.\Ltens H^1(\underline K_v,A)[-2]&\hfl{+1}{}\cr
R\G_{c}(\underline U,({\cal A}_{|U,p^.}))&\hfl{}{}&Sel_{p^.}(A/\underline K)&\hfl{}{}&\mathop\prod\limits_{v\in |Z|}\Z/p^.\Ltens A(\underline K_v)[-1]&\hfl{+1}{}}$$

\end{prop} Proof. $(i)$ The first and second distinguished triangles follow immediately from Rem. \ref{remselmer} $(ii)$ and Prop. \ref{selindepU} $(ii)$. The third one follows from the second one since  \' etale cohomology over $Spec(K)=\limp_{U\in \cal U} U$ can be computed as a limit by \cite{SGA4}, VII, 5.7 (note that the cohomology of $\cal A$ can be computed with respect to the \' etale topology since it is a smooth algebraic group).

$(ii)$ The first distinguished triangle is deduced from the first one of $(i)$ by applying $\Z/p^.\otimes^L(-)[-1]$ (observe that this functor commutes to $R\G(\underline C,-)$ and use (\ref{sp})). To get the second one, we remark that the commutative diagram $$\diagram{R\G(\underline C,\cal A)&\hfl{sp_Z}{}&R\G(\underline Z,\cal A_{|Z})\cr \vflcourte{d}{\parallel}&&\vflcourte{d'}{}\cr
R\G(\underline C,\cal A)&\hfl{sp^1_Z}{}&\tau_{\ge 1}R\G(\underline Z,\cal A_{|Z})\cr
\vflupcourte{c}{\parallel}&&\vflupcourte{c'}{}\cr
R\G(\underline C,\cal A)&\hfl{loc_C^1}{}& \mathop\prod\limits_{v\in |Z|}\tau_{\ge 1}R\G(\underline C_v,\cal A_{|C_v})}$$
in the derived category of $_{\underline G}Mod(\Z)$ has an canonical counterpart in the derived category of
diagrams of $_{\underline G}Mod(\Z)$ of this form since it only uses truncation and functoriality with
respect to the base. It follows from this remark that there is a canonical meaning for the mapping fiber of horizontal arrows and that there are canonical morphisms between them. The mapping fiber of the middle one is naturally isomorphic to $Sel(A/\underline K)$ (use Prop. \ref{selindepU} $(i)$, $(ii)$). By Prop. \ref{selindepU}$(ii)$, $c'$ is an isomorphism and we thus have $$Sel(A/\underline K)\simeq MF(loc^1_C)\simeq MF(sp_Z^1)$$
Now $d$ being an isomorphism, we get a canonical distinguished triangle $MF(sp_Z)\rightarrow MF(sp^1_Z)\rightarrow MF(d')\rightarrow MF(sp_Z)[1]$, ie. $$R\G^Z(\underline C,\cal A)\rightarrow
Sel(A/\underline K)\rightarrow H^0(Z,\cal A_{|Z})\mathop\rightarrow\limits^{+1}$$
The second distinguished triangle follows by applying $\Z/p^.\otimes^L(-)[-1]$.

$(iii)$  The first distinguished triangle is deduced from the second of $(i)$ by applying $\Z/p^.\otimes^L(-)[-1]$.  For the second one, we observe that the commutative diagram $$\diagram{R\G(\underline U,\cal A_{|U})&\hfl{loc^1_U}{}&
\mathop\prod\limits_{v\in |Z|}\tau_{\ge 1}R\G(\underline U_v,\cal A_{|U_v})\cr
\vflupcourte{e}{\parallel}&&\vflupcourte{e'}{}\cr
R\G(\underline U,\cal A_{|U})&\hfl{loc_U}{}&
\mathop\prod\limits_{v\in |Z|}R\G(\underline U_v,\cal A_{|U_v})\cr}$$
has an obvious counterpart in the derived category of such diagrams since it only uses truncation and functoriality with respect to the base. This justifies working with functorial mapping fibers as in $(i)$ and $(ii)$. We get a canonical distinguished triangle $MF(loc_U)\rightarrow MF(loc^1_U)\rightarrow MF(e')\rightarrow
MF(loc_U)[1]$, ie. $$R\G_c(\underline U,\cal A_{|U})\rightarrow Sel(A/\underline K)\rightarrow \mathop\prod\limits_{v\in |Z|}H^0(\underline U_v,\cal A_{|U_v})\mathop\rightarrow\limits^{+1}$$
The second distinguished triangle follows by applying $\Z/p^.\otimes^L(-)[-1]$.

\begin{prop} \label{duality}
Let $\hat A/K$ denotes the abelian variety which is dual to $A/K$. There is a canonical duality isomorphism in $D^{[0,3]}({_{\underline G}{Mod(\Zp)}}^{N_2})$:
$$Sel_{p^.}(\hat A/\underline K)\simeq Sel_{p^.}(A/\underline K)^\vee[-3]$$
where $(-)^\vee$ means the exact functor $Hom_{\Zp}(-,\Qp/\Zp)$.
\end{prop}
Proof. That $Sel_{p^.}(\hat A/\underline K)$ has no cohomology outside $[0,3]$ follows for instance from the third distinguished triangle in Prop. \ref{compsel} $(i)$ since $K$ has strict cohomological dimension $\le 3$.
Let us come to the duality isomorphism. A sketch of proof is given in \cite{KT} for a weaker statement (cf \emph{loc. cit} 2.4) by combining local duality over $U_v$, $v\in |Z|$ together with global duality over $U$. Here we concentrate on the ungrateful task of checking that it is indeed possible to prove the result in the setting of derived categories of normic systems. For this purpose referring to a collection of pairings at each level $n$ and checking compatibilities is not sufficient. The main point is thus to build the required duality morphism in $D^b({_{\underline G}{Mod(\Zp)}}^{N_2})$. That it is an isomorphism, will follow from the compatibility of our construction with that of \cite{Mi2}, III.\petit

For the purpose of the proof, let $B/K:={\hat A}/K$, $\cal B/C$ its N\'eron model and $\Psi/C$ its component group. Recall (\cite{SGA7} or \cite{Mi2}) that the Poincar\'e biextension $P:\cal A_{|U}\otimes^L\cal B_{|U}\rightarrow \Gm_{|U}[1]$ (which is essentially given by definition of $\hat B/F$) induces the Grothendieck pairing
$G:z^{-1}\Phi \otimes z^{-1}\Psi\rightarrow \Q/\Z$ (essentially via $j_*\Gm[1]\rightarrow z_*\Z[1]\leftarrow z_*\Q/\Z$). We know that $P$ extends over $C$ to a (unique) biextension $\cal A^{\Phi'}\otimes_L\cal B^{\Psi'}\rightarrow \Gm[1]$ as long as $\Phi'\subset \Phi$ and $\Psi'\subset \Psi$ are orthogonal with respect to $G$. Taking $(\Phi',\Psi')=(0,\Psi)$ and then $(\Phi,0)$, we thus get a canonical commutative square (of sheaves on the small smooth site $C_{sm}$)
\begin{eqnarray}\label{poinc}\diagram{\cal A^0&\hfl{}{}&\cal Ext^1_{C_{sm}}(\cal B,\Gm)\cr \vflcourte{}{}&&\vflcourte{}{}\cr
\cal A&\hfl{}{}&\cal Ext^1_{C_{sm}}(\cal B^0,\Gm)}\end{eqnarray}

Our proof consists in three steps.  \petit

\emph{Step 1.} The compatible pair of arrows (\ref{poinc}) induces morphisms as follows: \begin{eqnarray}\label{morpoinc}[\cal A^0\rightarrow \cal A]\rightarrow R\cal Hom_{C_{sm}}([\cal B^0\rightarrow \cal B],\Gm[1])& \hbox{in $D^b(C_{sm}^{[1]})$}\\
\label{morpoincpk}R\G(\underline C,[\cal A^0\rightarrow \cal A])\rightarrow  R\G(\underline C,[\cal B^0\rightarrow \cal B])^\vee[-2]& \hbox{in $D^b({_{\underline G}{Mod(\Z)}^{[1]}})$}\end{eqnarray}

Let us indicate how to get these. The arrows of (\ref{poinc}) give rise to
$$\begin{array}{lcr}&[\cal A^0\rightarrow \cal A]\rightarrow \cal Ext^1_{C_{sm}}([\cal B^0\rightarrow \cal B],\Gm)& \hbox{in $Mod(C_{sm}^{[1]},\Z)$}\\
\hbox{and thus to} &[\cal A^0\rightarrow \cal A]\rightarrow (\tau_{\ge 1}R\cal Hom_{C_{sm}}([\cal B^0\rightarrow \cal B],\Gm))[1]&\hbox{in $D^b(C_{sm})$}\end{array}$$ and (\ref{morpoinc}) follows since $\cal Hom_{C_{sm}}([\cal B^0\rightarrow \cal B],\Gm)=0$ by  \cite{SGA7}, 
VIII, 3.2 since $\Gm\hookrightarrow j_*j^*\Gm$. Applying $R\G(\underline C,-)$ we get $$R\G(\underline C,[\cal A^0\rightarrow \cal A])\rightarrow   R\underline{Hom}_{C_{sm}}([\cal B^0\rightarrow \cal B],\Gm)$$ where $\underline{Hom}_{C_{sm}}([M_1\rightarrow M_2],M_3)$ means the object of ${_{\underline G}{Mod(\Z)}^{[1]}}$ whose $n$th component is $[Hom_{C_{n,sm}}(M_{2|C_n},M_{3|C_n})\rightarrow Hom_{C_{n,sm}}(M_{1|C_n},M_{3|C_n})]$ with $G_n$ acting by conjugation and where going up (resp. down) transition maps are given by restriction (resp. trace). We note that there is a natural morphism
$$\underline{Hom}_{C_{sm}}([M_1\rightarrow M_2],M_3)\rightarrow Hom(\G(\underline C,[M_1\rightarrow M_2]),\Gamma(C,M_3))$$ which is bifunctorial with respect to $([M_1\rightarrow M_2],M_3)$ where for $([N_1\rightarrow N_2],N_3)$ in $\smash{{_{\underline G}{Mod(\Z)}^{[1]}}}\times Mod(\Z)$, $Hom([N_1\rightarrow N_2],N_3)$ means the normic system whose $n$th component is $[Hom(N_{2,n},N_3)\rightarrow Hom(N_{1,n},N_3)]$ and where going down (resp. up) transition maps are deduced from the going up (resp. down) transition maps of $N_2$ and $N_3$. This natural transformation gives rise to a map
$$R\underline{Hom}_{C_{sm}}([M_1\rightarrow M_2],M_3)\rightarrow RHom(\G(\underline C,[M_1\rightarrow M_2]),\Gamma(C,M_3))$$ in $D^-(\smash{{_{\underline G}{Mod(\Z)}^{[1]}}})$ as long as $R\G(C,M_3)$ is bounded. Applying this to $([\cal B^0\rightarrow \cal B],\Gm)$ and using the trace morphism $R\G(C,\Gm)\rightarrow \Q/\Z[-2]$ finally gives (\ref{morpoincpk}) as claimed.

It is worth noting that the pairings on cohomology resulting from (\ref{morpoincpk}) are compatible with those of \cite{Mi2}, III. To state this fact precisely, let us fix any closed point $v\in |C|$ as well as an open $U\subset C$ over which $\cal A$ has good reduction. Then it follows from the definitions that one has a natural commutative ``diagram''

$$\diagram
{H^q(Z_{n,v},\Phi_{|Z_{n,v}})&\times & H^{1-q}(Z_{n,v},\Psi_{|Z_{n,v}})&\hflcourte{1}{}&H^1(Z_{n,v},\Q/\Z)&&\cr
\vflupcourte{}{}&&\vflupcourte{}{}&&\vflupcourte{\wr}{}&&\cr
H^q(C_{n,v},\cal A_{|C_{n,v}})&\times & H^{1-q}(C_{n,v},\cal B_{|C_{n,v}})&\hflcourte{2}{}&H^2(C_{n,v},j_{v,*}\Gm)&&\cr
\vflcourte{}{}&&\vflcourte{}{}&&\vflcourte{\wr}{}&&\cr
H^q(U_{n,v},A_{|U_{n,v}})&\times & H^{1-q}(U_{n,v},B_{|U_{n,v}})&\hflcourte{3}{}&H^2(U_{n,v},\Gm)&&\cr
\vflupcourte{}{}&&\vflcourte{}{}&&\vflcourte{\wr}{}&&\cr
H^q(C_{n,v},\cal A_{|C_{n,v}})&\times & H^{2-q}_{Z_{n,v}}(C_{n,v},\cal B^0_{|C_{n,v}})&\hflcourte{4}{}&H^3_{Z_{n,v}}(C_{n,v},\Gm)&&\cr
\vflupcourte{}{}&&\vflcourte{}{}&&\vflcourte{\wr}{}&&\cr
H^q(C_n,\cal A_{|C_n})&\times &H^{2-q}(C_n,\cal B^0_{|C_n})&\hflcourte{5}{}& H^3(C_n,\Gm)&\hflcourte{Tr}{\sim}&H^3(C,\Gm)&\simeq & \Q/\Z\cr
\vflupcourte{}{}&&\vflcourte{}{}&&\vflupcourte{\wr}{}&&\cr
H^q_c(U_n,\cal A_{|U_n,p^k})&\times & H^{3-q}(U_n,\cal B_{|U_n,p^k}^0)&\hflcourte{6}{}&H^3_c(U_n,\Gm)}$$

where:

\noindent - pairings $1$ and $2$ are the ones occurring during the proof of \cite{Mi2}, III  Thm. 7.11 and right vertical isomorphism is by $j_{v,*}\Gm[1]\rightarrow z_{v,*}\Z[1]\leftarrow z_{v,*}\Q/\Z$.

\noindent - pairing $3$ (resp. $4$, $6$) pairing is the one occurring in \cite{Mi2}, III 7.8 (resp. 7.13, resp. 8.2).


\noindent - pairing $5$ is defined using the Yoneda pairing for Ext's and the vanishing of $\cal Hom_{C_n}(\cal B^0,\Gm)$ (or alternatively, the vanishing of $\cal Hom_{C_n}(\cal A,\Gm)$). By the restriction-corestriction formula, it coincides with the one induced by (\ref{morpoincpk});

It can be checked that similar compatibilities hold in the derived category. The proof however is not so easy to write down. We will thus avoid it and use the following trick instead. The morphism $(\ref{poinc})$ induces a map $$\Z/p^k\Ltens\cal A [-1]\rightarrow R\cal Hom_{C_{sm}}(\Z/p^k\Ltens\cal B^0,\Gm)$$
Now it it straightforward to check that this morphism induces compatible pairings $1',\dots,6'$ just as $1,\dots, 6$ above with $\cal A$ and $\cal B^0$ respectively replaced by  $\Z/p^k\otimes^L\cal A [-1]$ and $\Z/p^k\otimes^L\cal B^0$. Moreover, the resulting pairing
$$\diagram
{H^q(C_n,\Z/p^k\Ltens\cal A_{|C_n}[-1])&\times &H^{2-q}(C_n,\Z/p^k\Ltens\cal B^0_{|C_n})&\hflcourte{5'}{}& H^3(C_n,\Gm)&\hflcourte{Tr}{\sim}&H^3(C,\Gm)&\simeq & \Q/\Z}$$
coincides with the one deduced from (\ref{morpoincpk}) by applying $\Z/p^k\otimes^L(-)$.

\petit

\emph{Step 2.} The morphism (\ref{morpoinc}) induces a canonical morphism of distinguished triangles in $D^b({_{\underline G}Mod(\Z)})$
$$\diagram{Sel(A/\underline K)&\hfl{}{}&R\G(\underline C,\cal A)&\hfl{}{}&H^1(\underline Z,\Phi_{|Z})[-1]&\hfl{+1}{}\cr
\vflcourte{a}{}&&\vflcourte{b}{}&&\vflcourte{c}{}\cr
Sel(B/\underline K)^\vee[-2]&\hfl{}{}&R\G(\underline C,\cal B^0)^\vee[-2]&\hfl{}{}&H^0(\underline Z,\Psi_{|Z})^\vee[-1]&\hfl{+1}{}}$$

\petit

Let us explain this. Since it can be expressed in terms of truncation, functorial cones and natural transformation of functors in the variable $[M\mathop\rightarrow\limits N]\in Mod(C_{sm},\Z)^{[1]}$, we observe that the following diagram of $D({_{\underline G}Mod(\Z)})$:  $$\diagram{R\G(\underline C,N)&\hfl{\beta}{}&\tau_{\ge 1}Cone(R\G(\underline C,[M\rightarrow N]))\cr
\vflcourte{}{}&&\vflcourte{}{}\cr
R\G(\underline C,N)&\hfl{\beta'}{}&\tau_{\ge 1}Cone(R\G(\underline Z,[z^{-1}M\rightarrow z^{-1}N]))\cr
\vflupcourte{}{}&&\vflupcourte{}{}\cr
R\G(\underline C,N)&\hfl{loc_C^1}{}&\mathop\prod\limits_{v\in |Z|}\tau_{\ge 1}R\G(\underline C_v,N)}$$
has an obvious counterpart in the corresponding derived category of diagrams of this form.
Applying this to $[M\rightarrow N]=[\cal A^0\rightarrow \cal A]$ and forming  mapping fibers of the horizontal arrows, we get canonical isomorphisms (using arguments already explained in the proof of Prop. \ref{compsel}): $$MF(\beta)\mathop\leftarrow\limits^\sim MF(\beta')\mathop\rightarrow\limits^\sim MF(loc_C^1)\simeq Sel(A/\underline K)$$

Consider the following diagram as a functor whose argument $[\underline M\mathop\rightarrow\limits^f \underline N]$ varies in the category of complexes $C({_{\underline G}Mod(\Z)^{[1]}})$:

$$\diagram{\tau_{\le 1}MF(f)&\simeq& \tau_{\le 1}MF(f) \cr \vfl{}{}&& \vfl{\underline \alpha}{}\cr
MF(f)&\hfl{}{}&\underline M&\hfl{f}{}&\underline N&\hfl{}{}&Cone(f)\cr
\vfl{}{}&&\vfl{}{}&&\vfl{}{}&&\vfl{}{}\cr
\tau_{\ge 2}MF(f)&\hfl{}{}&MF(\underline \beta)&\hfl{}{}&\underline N&\hfl{\underline \beta}{}&\tau_{\ge 1}Cone(f)}$$

We remark the following:

\noindent - since the left vertical composed morphism is zero, so is the middle vertical one. Whence a canonical morphism $Cone(\underline \alpha)\rightarrow MF(\underline \beta)$. As shown by the diagram, this is automatically a quasi-isomorphism.

\noindent - the functors $Cone(\underline \alpha)$ and $MF(\underline \beta)$ are exchanged by the exact contravariant autofunctor $(-)^\vee=Hom_\Z(-,\Q/\Z)$ of $C({_{\underline G}Mod(\Z)^{[1]}})$. Indeed: $$\begin{array}{rcl}(Cone(\alpha)[M\rightarrow N])^\vee&\simeq & (Cone([\tau_{\le 1} MF(f)\rightarrow M]))^\vee\\
&\simeq & MF([M^\vee\rightarrow \tau_{\ge 1} Cone(f^\vee)])\\
&\simeq & MF(\beta)[N^\vee\rightarrow M^\vee]\end{array}$$


Putting everything together we get in $D({_{\underline G}Mod(\Z)})$: $$\begin{array}{rclc}Sel(A/\underline K)&\simeq & MF(\beta)[\cal A^0\rightarrow \cal A]&\\
& \simeq & MF(\underline \beta)(R\G(\underline C,[\cal A^0\rightarrow \cal A]))&\\
&\rightarrow&MF(\underline \beta)(R\G(\underline C,[\cal B^0\rightarrow \cal B])^\vee[-2]) &\hbox{(cf. (\ref{morpoincpk}))}\\
&\simeq &(Cone(\underline \alpha)(R\G(\underline C,[\cal B^0\rightarrow \cal B])))^\vee[-2]&\\
&\simeq &(MF(\underline \beta)(R\G(\underline C,[\cal B^0\rightarrow \cal B])))^\vee[-2]&\\
&\simeq &(MF(\beta)[\cal B^0\rightarrow \cal B])^\vee[-2]&\\
&\simeq &Sel(B/\underline K)^\vee[-2]&\end{array}$$

\emph{Step 3.} The morphisms $a$, $b$, $c$ of Step 2 become invertible when applied  $\Z/p.\otimes^L(-)$.\petit

It is of course sufficient to prove this for $b$ and $c$. By \cite{Mi2}, III, 7.11, $c$ itself is an isomorphism. Now the statement about $b$ is equivalent to the horizontal arrow $b'$ below being an isomorphism:

$$\diagram{\vflcourte{}{}&&\vflcourte{}{}\cr
H^q_c(U_n,\Z/p^k\Ltens\cal A_{|U_n}[-1])&\hfl{}{}&H^{2-q}(U_n,\Z/p^k\Ltens\cal B^0_{|U_n})^\vee\cr
\vflcourte{}{}&&\vflcourte{}{}\cr
H^q(C_n,\Z/p^k\Ltens\cal A_{|C_n}[-1])&\hfl{b'}{}&H^{2-q}(C_n,\Z/p^k\Ltens\cal B^0_{|C_n})^\vee\cr
\vflcourte{}{}&&\vflcourte{}{}\cr
\mathop\prod\limits_{v\in |Z_n|} H^q(C_{n,v},\Z/p^k\Ltens\cal A_{|C_{n,v}}[-1])&\hfl{}{}&
\mathop\prod\limits_{v\in |Z_n|} H^{2-q}_{Z_{n,v}}(C_{n,v},\Z/p^k\Ltens\cal B^0_{|C_{n,v}})^\vee
}$$

This diagram commutes thanks to the compatibilities explained in Step 1 and its columns are long exact sequences as follows easily from the definition of compactly supported cohomology together with the complete excision property
$$\mathop\prod\limits_{v\in |Z_n|} H^{2-q}_{Z_{n,v}}(C_{n,v},\Z/p^k\Ltens\cal B^0_{|C_{n,v}})\simeq H^{2-q}_{Z_{n}}(C_{n},\Z/p^k\Ltens\cal B^0_{|C_{n}}),$$
obtained by using \cite{Mi1}, III, 1.28 and \cite{Mi2}, III, 7.14. The first (resp. third) horizontal morphism is invertible thanks to \cite{Mi2}, III, 8.2 (resp. \cite{Mi2}, III 7.13) and the result follows.
\begin{flushright}$\square$\end{flushright}

\begin{rem} The proof does not require any assumption about the reduction of $A$ at the points of $Z$.  \end{rem}

\begin{cor} \label{selcoh} Consider the usual Selmer group $$Sel^1_{p^k}(A/K_n):=Ker \left(H^1(K_n,A_{|K_n,p^k})\rightarrow \mathop\oplus\limits_{v\in |C|}H^1(K_{n,v}A_{|K_{n,v}}))\right)$$ and let $Sel^1_{p^.}(A/\underline F)$ denote the object of $Mod({_{\underline G}Mod(\Z)^{N_2}})$ encapsulating the whole collection (together with the various natural transition morphisms) for varying $n,k$. The cohomology of the complex $Sel_{p^.}(A/\underline K)$ vanishes outside $[0,3]$ and is described as follows in this range:
$$\begin{array}{rcl}H^0&=&A_{p^.}(\underline K)\cr
H^1&=&Sel^1_{p^.}(A/\underline K)\\
H^2&\simeq &Sel^1_{p^.}(\hat A/\underline K)^\vee\\
H^3&\simeq &\hat A_{p^.}(\underline K)^\vee\end{array}$$
\end{cor}

Proof. The description of $H^0$ and $H^1$ follows easily from the following distinguished triangle (Prop. \ref{compsel} $(i)$):
$$\diagram{Sel_{p^.}(A/\underline K)&\hfl{}{}R\G(\underline K,A_{p^.})&\hfl{}{}&\mathop\oplus\limits_{v\in |C|}\Z/p^.\Ltens H^1(\underline{U}_v,A_{|U_v})[-2]&\hfl{+1}{}\cr}$$
whose associated long exact sequence of cohomology reads as the collection of the following ones for varying $n$, $k$:  $$\diagram{0&\hfl{}{}&H^0(Sel_{p^k}(A/K_n))&\hfl{}{}&H^0(K_n,A_{p^k})&\hfl{}{} &0\cr
&\hfl{}{}&H^1(Sel_{p^k}(A/K_n))&\hfl{}{}&H^1(K_n,A_{p^k})&\hfl{}{}&\mathop\oplus\limits_{v\in |C|}Tor_1(\Z/p^k,H^1(K_{n,v},A_{|K_{n,v}}))&\dots}$$
(note that for any abelian group $M$, $Tor_1(\Z/p^k,M)$ is identified with the kernel of $p^k$ on $M$. It it easy to check that modulo this identification, the right arrow in the second exact sequence is induced from the localization map occurring in the definition of $\smash{Sel^1_{p^k}}(A/K_n)$).

The description of $H^2$ and $H^3$ follows by duality, thanks to Prop. \ref{compsel} and Prop. \ref{duality}.
\begin{flushright}$\square$\end{flushright}

\begin{cor} In $D^b({_{\underline G}{Mod(\Zp)}})$, one has the duality isomorphism
$$Sel_{T_p}(\hat A/\underline K)\simeq Sel_{p^\infty}(A/\underline K)^\vee[-3]$$  and those complexes are concentrated in $[1,3]$.
\end{cor}
Proof. Just apply $R\limp_k$ to Prop. \ref{duality} and notice that $H^0$ vanishes because for $n$ fixed,  $\hat A_{p^k}(K_n)$ is bounded
independently of $k$.
\begin{flushright}$\square$\end{flushright}

\begin{cor} \label{arcoh} $R\G^Z(\underline C,\cal A_{p^.})$ is concentrated in $[0,3]$ and its cohomology objects fit into compatible exact sequences and isomorphisms as follows:
$$\begin{array}{l}H^0\hookrightarrow A_{p^.}(\underline K)\rightarrow \cal A_{p^.}(\underline Z)
\rightarrow H^1\rightarrow Sel_{p^.}(A/\underline K)\twoheadrightarrow \cal A(\underline Z)/p^.\\ H^2\simeq
Sel_{p^.}(\hat A/\underline K)^\vee\\
H^3\simeq \hat A_{p^.}(\underline K)^\vee\end{array}$$
 \end{cor} \medskip

\begin{rem}\label{SelGar} Even though we won't use it, let us point out the relation to the complexes of \cite{KT}.

-  The complexes $Sel_{p^\infty}(A/K)=\limi_k Sel_{p^k}(A/K)$ and $R\G_{ar}\{p\}$ have the same cohomology (compare \emph{loc. cit.} Lem. 2.4 and our Cor. \ref{selcoh}).

- The cohomology of the complexes $\limi_k R\G^Z(C,\cal A_{p^k})$ and $R\G_{ar,V}\{p\}$ fits into similar exact sequences (compare \cite{KT}, 2.5.2, 2.5.3 and our Prop. \ref{compsel} (ii)) for $V=(V_v)_{v\in Z}$ and  $V_v=Ker(A(K)\rightarrow \cal A(Z_v))$.

- In fact, one may easily form isomorphisms $$Sel_{p^\infty}(A/K)\simeq R\G_{ar}\{p\}$$ and $$\limi_k R\G^Z(C,\cal A_{p^k})\simeq R\G_{ar,V}\{p\}.$$
\end{rem}

\grand

\subsection{Crystalline syntomic complexes} ~~ \label{reminderTV} \\

In this section we assume that $A/K$ is semistable and we recall the comparison result of \cite{TV} and deduce from it the fundamental distinguished triangles (\ref{tdfundnorm}), announced in the introduction, relating (normic) flat cohomology vanishing at $Z$ to (normic) crystalline cohomology.

\para We keep the notations of the previous section ($A/K$, ${\mathcal A}/C$, $U$, $Z$, $\epsilon: C_{FL}\to C_{et}$) and we let furthermore $C^\sharp$ denote the log-scheme whose underlying scheme is $C$ and whose log-structure is induced by $Z$. In the chapter 8 of \cite{TV} is associated to the semistable abelian variety $A/K$ a log Dieudonn\'e crystal $(D,F,V)$. Recall that such object is a crystal of locally free $\mathcal O$-modules of finite rank in the crystalline small \'etale ringed topos $((C^\sharp/\Zp)_{crys,et},{\mathcal O})$ endowed with two operators $F:\sigma^*D\to D$, $V:D\to \sigma^*D$ (where $\sigma$ denotes the endomorphism of $(C^\sharp/\Zp)_{crys,et}$ induced by the absolute Frobenius of $C^\sharp$) such that $FV=VF=p$. This log-crystal extends the covariant Dieudonn\'e crystal $(D_{U}({\mathcal A}|_U),F,V)$ of \cite{BBM} and is endowed  with a canonical epimorphism
$$\pi:D\to Lie({\mathcal A}),$$
where the Lie algebra $Lie({\mathcal A})$ is seen as a crystalline sheaf (\cite{TV}, 5.35 (i)).\\

We set $Fil^1D:=Ker(\pi)$ and denote ${\bf 1}:Fil^1 D\to D$ the canonical injection.

\para Let $u:(C^\sharp/\Zp)_{crys,et}\to C_{et}$ denote the canonical projection of the crystalline topos on $C_{et}$. Following \cite{TV}, Prop. 5.49, we recall the construction of an operator
$$\phi:Ru_*Fil^1D(-Z)\to Ru_*D(-Z),$$
where the twist by $-Z$ of a crystalline sheaf is defined in \cite{TV}, Def. 4.42.

Consider $(U_{[.]}^\sharp, Y_{[.]}^\sharp,\iota_{[.]},\tilde \sigma)$ where  $U_{[.]}^\sharp$ is a semisimplicial fine log-scheme above $C^\sharp$, $Y_{[.]}^\sharp$ is a semisimplicial $p$-adic formal log scheme with finite local $p$-bases over $Spf(\Zp)$ (in the sense of \cite{TV}, 4.1 (v)), endowed with a lifting of Frobenius $\tilde \sigma$ and $\iota:U_{[.]}^\sharp\to Y_{[.]}^\sharp$ is a closed immersion. We denote $T_{[.],n}^\sharp$ the log pd envelope of $U_{[.]}^\sharp$ into the reduction $Y_{[.],n}^\sharp$ of $Y_{[.]}^\sharp$ modulo $p^n$ and, for a crystalline sheaf $E$, $E_{T_{[.],n}^\sharp}$ the realization of the sheaf at the pd-thickening $(U_{[.]}^\sharp,T_{[.],n}^\sharp)$. It is proved (\cite{TV}, Lem. 4.9) that we can always chose these data such that $(U_{[.]}^\sharp,T_{[.],n}^\sharp)$ is a hypercovering in the topos $(C^\sharp/\Zp)_{crys,et}$. In \cite{TV}, 5.49 it is proved that there exists a unique collection $(\phi_n)$ of $\tilde \sigma$-semi-linear morphisms of $\mathcal O_{T_{[.],n}}$-modules rendering the following diagram commutative for all $n$:
$$\xymatrix{{{Fil{}^1D(-Z)_{T_{[.],n+1}^\sharp}}/p^n}\ar[r]^-{{\phi}_n}\ar[d]^1&{
{D(-Z)_{T_{[.],n}^\sharp}}}\ar[d]^p\\
{{D(-Z)_{T_{[.],n}^\sharp}}}\ar[r]^-{Fr}&(\tilde F)_*{{D(-Z)_{T_{[.]}^\sharp}}}}$$
(here $Fr$ is the morphism built up from the operator $F$ of $D$ together with the lifting of frobenius $\tilde \sigma$ of $T_n$).
Finally, we set $$\phi:=R\ilim{} Rf_{T_{[.],n},*}{\phi}_n,$$
where $f_{T_{[.],n}}:T_{[.],n,et}\to C_{et}$. By cohomological descent in the crystalline topos, this yields a morphism (\cite{TV}, Prop. 5.55, (ii) and Rem. 5.56)
$$\phi:Ru_*Fil^1D(-Z)\to Ru_*D(-Z).$$

\para \label{paradt} To simplify the notations, we set:
\begin{enumerate}
\item $N:=R\epsilon^Z_* R\ilim \mathcal{A}_{p^n}$,
\item $I:=Ru_*Fil^1D(-Z)$,
\item $P:=Ru_*D(-Z)$,
\item $L:=Lie(\mathcal{A})(-Z)$.
\end{enumerate}
All $N,I,P,L$ are considered as objects of the derived category $D(C_{et},\Zp)$.

The main theorem of \cite{TV}  implies (note that $\epsilon^Z_*$ as defined in \ref{defvanishingcoh} coincides with the functor $\epsilon_* \underline \Gamma^Z$ used in \emph{loc. cit.}):

\begin{thm}\label{key-dt}(\cite{TV}, Cor. 9.17)
In $D^p(C_{et},\Zp)$ there are canonical distinguished triangles as follows:

$$N\to I\buildrel{1-\phi}\over\to P\buildrel{+1}\over\to$$
$$I\buildrel{1}\over\to P\to L\buildrel{+1}\over\to$$
\end{thm}
\begin{flushright}
$\square$
\end{flushright}

\para  \label{paradtinfty}
For $X\in{N,I,P,L}$, we denote $\underbar X:=(X_n)$ the object $R\Gamma(\underline{C},X)\in D^b({_{\underline G}(\Zp\mod)})$ (so that $X_n=R\G(C_n,X_{|C_n})$, see Sect. \ref{sectexamples} for explanations). Given a Galois subextension $K'/K$ of $K_\infty/K$, we denote $X_{K'}$ (resp. $X_{K'}^*$) the derived projective limit of the projective system of complexes indexed by the finite subextensions of $K'/K$ (eg. $X_{K_n}=X_n$, $X_{K_n}^*=X_n^*$, $X_{K_\infty}:=R\ilim X_n$,  $X_{K_\infty}^*:=R\ilim RHom_{\Zp}(X_n,\Zp)$). We have by Thm.  \ref{key-dt} two pairs of distinguished triangles in $D(\La)$:

\begin{enumerate}\label{dtinfty}
\item $N_{K_\infty}\to I_{K_\infty}\buildrel{1-\phi}\over\to P_{K_\infty}\buildrel{+1}\over\to$
\item $I_{K_\infty}\buildrel{1}\over\to P_{K_\infty}\to L_{K_\infty}\buildrel{+1}\over\to$
\end{enumerate}
and
\begin{enumerate}\label{dtinfty*}
\item $P^*_{K_\infty}\buildrel{1^*-\phi^*}\over\to I^*_{K_\infty}\to N^*_{K_\infty}\buildrel{+1}\over\to$
\item $L^*_{K_\infty}\to P^*_{K_\infty}\buildrel{1}\over\to I^*_{K_\infty}\buildrel{+1}\over\to$
\end{enumerate}

\begin{rem}\label{remTV} These distinguished triangles are obtained from Thm. \ref{key-dt} by applying $R\G(\underline C,-)$ then $R\ilim$. Note that, even though \cite{KT} might give distinguished triangles in the derived category of modules at each level $n$ and transition mophisms between them, this would formally not be enough to apply  $R\ilim$ (whose source needs to be the derived category of projective systems rather than the category of projective systems of the derived category) and get the desired distinguished triangles in $D(\La)$. One could of course easily give a meaning to $R\ilim$ of each term of the triangle using explicit complexes, but this is far from clear for the morphisms, especially for the morphism $N_n\to I_n$.
\end{rem}

In fact, the two cohomology theories $X_{K_\infty}$ and $X_{K_\infty}^*$ are related as follows:

\begin{prop}\label{La-duality} Let $X\in{N,I,P,L}$. Then $X_{K_\infty}$ and $X_{K_\infty}^*$ are in $D^p(\La)$. Moreover:

$(i)$ If $\underline \La=(\Z[G_n])$ denotes the canonical normic system then we have in $D({_{\underline G}(\Zp\hbox{-}mod}))$ isomorphisms as follows:
$$\underline \L\Ltens_\L X_{K_\infty}\simeq \underline X \hspace{1cm} \hbox{and} \hspace{1cm} \underline \L\Ltens_\L X^*_{K_\infty}\simeq \underline X^*$$
(so that, in particular, we have in $D(_{G_n}(\Z_p\hbox{-}mod))$:
$$\Z_p[G_n]\Ltens_\L X_{K_\infty}\simeq X_n \hspace{1cm} \hbox{and} \hspace{1cm} \Z_p[G_n]\Ltens_\L X^*_{K_\infty}\simeq  X_n^*)$$
$(ii)$ We have an isomorphism in $D^p(\La)$:
$$RHom_{\La}(X^*_{K_\infty},\La)\simeq X_{K_\infty}.$$
\end{prop}

Proof. Since $C$ is proper and smooth, well known finiteness results imply that $\underline L$ and $\underline P$ satisfy the finiteness condition of Thm. \ref{main1}. The perfectness statement follows directly from Thm. \ref{main1} for $X=L$ or $P$. It then follows from the distinguished triangle $(ii)$ (resp. $(i)$) for $X=I$ (resp. $N$).
The remaining statements follow from Thm. \ref{main1} as well (using \cite{Va} 2.10.3 and 3.1).\begin{flushright}$\square$\end{flushright}

\section{The Main conjecture}

\subsection{The constant $\Zp$-extension} \hbox{  } ~~ \\

Recall that $K_{ar}:=Kk_\infty$ denotes the constant $\Zp$-extension of $K$. During this paragraph, we index by $k_n$ the $n$-th layer of
$k_\infty/k$.

\begin{prop} \label{propbc} Let $\underline W:=(W(k_n))\in {_{\underline G}(\Zp\mod)}$ denote the normic system of $\Zp$-modules formed by
the Witt vectors along $k_\infty/k$. Also, let $\underline W\otimes_\Zp(-):\Zp\mod\rightarrow {_{\underline \G}(\Zp\mod)}$ denote the obvious
functor. Then:

$(i)$ For $X\in \{I,P,L\}$, there is a canonical isomorphism in $D^b({_{\underline G}(\Zp\mod)})$: $$\underline W\Ltens X_K\simeq \underline
X$$

$(ii)$ The morphism   ${\bf 1}:I_n\rightarrow P_n$ (resp. $\phi:I_n\rightarrow  P_n$)  is $W(k_n)$-linear (resp. $W(k_n)$-semi-linear with respect to
the action of the absolute frobenius $\sigma:W(k_n)\rightarrow W(k_n)$). In particular, the morphism ${\bf 1}-\phi:\underline I\rightarrow \underline
P$ can be described at level $n$ as $\lambda\otimes x_0\mapsto \lambda\otimes {\bf 1} x_0-\sigma\lambda\otimes \phi x_0$.
\end{prop}
Proof. $(i)$ follows from the Zariski log-crystalline base change theorem (note that the base change map is well defined at the level complexes and
thus makes sense in $D^b{_{\underline \G}(\Zp\mod)}$).

$(ii)$ Linearity of ${\bf 1}$ and semi-linearity of $\phi$ are immediate from the definitions. The second statement
folllows.
\begin{flushright}$\square$\end{flushright}

\begin{cor} \label{torsionar} Both $N_{K_{ar}}$ and $N^*_{K_{ar}}$  are torsion over $\La(\Gamma)$, ie. $$Q(\Gamma)\Ltens_{\La(\Gamma)}N_{K_{ar}}=0
\hspace{1.5cm}\hbox{and}\hspace{1.5cm}Q(\Gamma)\Ltens_{\La(\Gamma)}N^*_{K_{ar}}=0$$
\end{cor}
Proof. Since by Prop. \ref{La-duality}, we have $N_{K_{ar}}\simeq RHom_{\La}(N^*_{K_{ar}},\La)$, it suffices to prove that $N_{K_{ar}}$ is torsion. Applying
$\Qp\otimes_\Zp(-)$ to the distinguished triangles of \ref{paradt} gives rise to the distinguished triangle:
$$\Qp\otimes_\Zp N_{K_{ar}}\rightarrow \Qp\otimes_\Zp I_{K_{ar}}\buildrel{1-\phi}\over\rightarrow \Qp\otimes_\Zp P_{K_{ar}}\rightarrow \Qp\otimes_\Zp N_{K_{ar}}[1]$$
which can be rewritten as
$$\Qp\otimes_\Zp N_{K_{ar}}\rightarrow W_\infty\Ltens_\Zp (\Qp\otimes_\Zp I_{K})\rightarrow
W_\infty\Ltens_\Zp(\Qp\otimes_\Zp P_{K})\rightarrow \Qp\otimes_\Zp N_{K_{ar}}[1]$$ Now $\Qp\otimes_\Zp L_{K_{ar}}=0$, since $Lie(D)$ is an $\mathbb
F_p$-vector space and ${\bf 1}: \Qp\otimes_\Zp I_{K_{ar}}\rightarrow \Qp\otimes_\Zp P_{K_{ar}}$ thus has an inverse. Since moreover $\Qp$ and
$W_\infty:=\limp W(k_n)$ are $\Zp$-flat, the long exact sequence of cohomology of the previous triangle can be written
$$\label{sel}\dots\to \Qp\otimes_\Zp H^q(N_{K_{ar}})\to W_\infty\Ltens_\Zp(\Qp\otimes_\Zp
H^q(P_{K}))\buildrel{id-\phi{\bf 1}^{-1}}\over\to W_\infty\Ltens_\Zp(\Qp\otimes_\Zp H^q(P_{K}))\to\dots$$ where
$(id-\phi{\bf 1}^{-1})(\lambda\otimes x)=\lambda \otimes x- \sigma\lambda\otimes \phi {\bf 1}^{-1}x$.

To end the proof, we need the following result from $\sigma$-linear algebra:

\begin{lem}  If $\psi$ is a linear endomorphism of a finite dimensional $\Qp$-vector space $M$ then:

$(i)$ The kernel of $id-\sigma\otimes \psi$ on $W(k_n)\otimes_\Zp M$ has $\Qp$-dimension bounded by the number of unit eigenvalues of $\psi$
(counted with multiplicities). In particular, its dimension is bounded independently of $n$.

$(ii)$ The $\Qp$-linear endomorphism $id-\sigma\otimes \psi$ of $(\limp W(k_n))\otimes_\Zp M)$ is injective.
\end{lem}
Proof. $(i)$ Set $L=Frac(W(\overline{\mathbb F_p}))$. Since $L^{\sigma}$ (the fixed points by $\sigma$) are reduced to $\Qp$ and
$W(k_n)\otimes_\Zp M\subset L\otimes_\Qp M$, \cite{EL2} 6.2 gives the result.

$(ii)$ Let $\G$ act on $\subset W(k_n)\otimes_\Zp M$ via $W(k_n)$. Thanks to $(i)$, we know that there is a open subgroup $\G_m\subset \G$ which fixes
the kernel of $id-\sigma\otimes \psi$ acting on $W(k_n)\otimes_\Zp M$ for all $n$. The result follows, since $$(W_\infty\otimes_\Zp M)^{\G_m}= (\limp
W(k_n)^{\G_m})\otimes_\Zp M=0.$$
\begin{flushright}$\square$\end{flushright}

Applying the lemma to $M=\Qp\otimes_\Zp H^q(P_{K})$ and $\psi=\phi{\bf 1}^{-1}$ shows that the first distinguished triangle of \ref{paradtinfty} reduces to short exact sequences $$0\rightarrow \Qp\otimes_\Zp
H^q(P_{K_{ar}})\buildrel{id-\psi}\over\rightarrow \Qp\otimes_\Zp H^q(P_{K_{ar}})\rightarrow \Qp\otimes_\Zp H^q(N_{K_{ar}})\rightarrow 0$$ and the result follows by applying $Q(\G)\otimes_{\La(\G)}(-)$
to those, since $Q(\G)\otimes_{\La(\G)} H^q(P_{K_{ar}})$ is a finite dimensional vector space over $Q(\G)$.
\begin{flushright}$\square$\end{flushright}

Following \cite{Ja}, we set $E^q(-):=Ext^q_{\La(\G)}(-,\La(\G))$. In loc. cit. these functors (and others) are used to classify $\La(\G)$-modules up to
isomorphism. Let us recall that, for a finitely generated $\La(\G)$-module they verify the following properties:

- $E^0(M)$ is free. It is zero if and only if $M$ is torsion.

- $E^1(M)$ is torsion and is (non canonically) pseudo-isomorphic to $M$. It has no non zero finite submodules if $M$ is torsion.

- $E^2(M)$ is finite.

\begin{cor} \label{descNar} $N_{K_{ar}}$ and $N_{K_{ar}}^*$ can be described as follows:

\noindent 1. $(i)$ $N_{K_{ar}}$ is concentrated in degrees $[1,3]$ and its cohomology $\L(\G)$-modules are described by the following exact sequence and isomorphisms
$$\begin{array}{l} H^1(N_{K_{ar}}) \hookrightarrow \limp_{n} Sel_{T_p}(A/Kk_n)\twoheadrightarrow \limp_{k,n} \prod_{v\in Z}\cal A(k_vk_n)/p^k\\
H^2(N_{K_{ar}})\simeq Hom_{\Zp}(\limi_{n} Sel_{p^\infty}(\hat A/Kk_n),\Qp/\Zp)\\
H^3(N_{K_{ar}})\simeq Hom_{\Zp}(\limi_{n} \hat A_{p^\infty}(Kk_n),\Qp/\Zp)\end{array}$$

$(ii)$ Moreover, one has the following isomorphisms and exact sequences of $\L(\G)$-modules:
$$\begin{array}{ccccc}&&H^1(N_{K_{ar}})&\simeq& E^1(H^0(N_{K_{ar}}^*))\\
E^2(H^0(N_{K_{ar}}^*))&\hookrightarrow & H^2(N_{K_{ar}})&\twoheadrightarrow &
E^1(H^{-1}(N_{K_{ar}}^*))\\
E^2(H^{-1}(N_{K_{ar}}^*))&\hookrightarrow & H^3(N_{K_{ar}})&\twoheadrightarrow&
E^1(H^{-2}(N_{K_{ar}}^*))\end{array}$$

\noindent 2. $(i)$ $N_{K_{ar}}^*$ is concentrated in degrees $[-2,0]$ and its cohomology $\L(\G)$-modules are described by the following exact sequence and isomorphism:
$$\begin{array}{l} H^{-2}(N_{K_{ar}}^*)\simeq \limp_{n} Sel_{T_p}(\hat A/Kk_n)\\
Hom_{\Zp}(\limi_{n}Sel_{p^\infty}(A/Kk_n),\Qp/\Zp)\hookrightarrow H^{-1}(N_{K_{ar}}^*) \rightarrow
Hom_{\Zp}(\limi_{n} \oplus_{v\in Z}\cal A_{p^\infty}(k_vk_n),\Qp/\Zp) \\
 \hspace{7cm}\rightarrow
Hom_{\Zp}(\limi_{n} A_{p^\infty}(Kk_n),\Qp/\Zp)\twoheadrightarrow H^0(N_{K_{ar}}^*)\end{array}$$

$(ii)$ Moreover, one has the following isomorphisms and exact sequences of $\L(\G)$-modules:
$$\begin{array}{ccccc}&&H^{-2}(N_{K_{ar}}^*)&\simeq& E^1(H^3(N_{K_{ar}}))\\
E^2(H^3(N_{K_{ar}}))&\hookrightarrow & H^{-1}(N_{K_{ar}}^*)&\twoheadrightarrow &E^1(H^2(N_{K_{ar}}))\\
E^2(H^2(N_{K_{ar}}))&\hookrightarrow & H^0(N_{K_{ar}}^*)&\twoheadrightarrow &E^1(H^1(N_{K_{ar}}))\end{array}$$
\end{cor}
Proof. The description of cohomology follows from Cor. \ref{arcoh} (recall that $\underline N=R\limp_k R\G^Z(\underline C,\cal A_{p^k})$ while $\underline N^*=RHom_\Zp(\underline N,\Zp)\simeq RHom_\Zp(\limi_kR\G^Z(\underline C,\cal A_{p^k}))$) using only that $A_{p^\infty}(Kk_n)$ and $\cal A_{p^\infty}(k_vk_n)$ are finite. The stated relations between $N_{K_{ar}}$ and $N_{K_{ar}}^*$ follow from Thm. \ref{La-duality} thanks to the spectral sequence
$$\begin{array}{lcccr}&Ext^p_{\La(\G)}(H^{-q}(N_{K_{ar}}^*),\La(\G))&\Rightarrow& Ext^{p+q}_{\La(\G)}(N_{K_{ar}}^*,\La(\G))& \\
\hbox{(resp.\hspace{1cm}}&Ext^p_{\La(\G)}(H^{-q}(N_{K_{ar}}),\La(\G))&\Rightarrow& Ext^{p+q}_{\La(\G)}(N_{K_{ar}},\La(\G))&\hbox{\hspace{.5cm})}\end{array}$$
which degenerates by Cor. \ref{torsionar}.
\begin{flushright}$\square$\end{flushright}

\begin{prop} \label{equivmu}The following conditions are equivalent:

$(i)$ $N_{K_{ar}}$ is finitely generated over $\Zp$ (that is for any $i$, $H^i(N_{K_{ar}})$ is a finitely generated $\Zp$-module).

$(ii)$ $Hom_\Zp(\limi_n Sel(\hat A_{p^\infty}/Kk_n),\Qp/\Zp)$ is finitely generated over $\Zp$.

$(iii)$ $N_{K_{ar}}^*$ is finitely generated over $\Zp$.

$(iv)$  $Hom_\Zp(\limi_n Sel(A_{p^\infty}/Kk_n),\Qp/\Zp)$ is finitely generated over $\Zp$.
\end{prop}

Proof. $(i)\Leftrightarrow (iii)$ follows from the previous corollary. Indeed if $M$ is any $\Lambda(\G)$-module then $E^i(M)$ is finitely generated over $\Zp$ as soon as $M$ is. The remaining equivalence are obvious, because the $p$-rank of torsion points remains bounded along the tower.


%
%
%
%
%
%
\begin{flushright}$\square$\end{flushright}

From now on, we will often need to make the following assumption: \medskip

\noindent{\bf ($\mu=0$ - up to iso)}: There exists an abelian variety $A'/K$ which is  isogenous to $A/K$ and verifies the equivalent conditions of Prop. \ref{equivmu}.


\begin{rem} $(i)$ \label{remmu} What we know about this assumption is the following:

- It holds for constant varieties except for supersingular abelian variety with non-invertible Hasse-Witt matrix (see \cite{OT}, Theorem 1.8).

- If it holds for $A/K$ and $K'/K$ is a finite unramified extension, then it also holds for $A/K'$.

$(ii)$ In characteristic $0$, it is generally expected that the analogue of ($\mu=0$-up to iso) always holds for elliptic curves (\cite{Gr}, p. 9).
\end{rem}

\subsection{The case of a $p$-adic Lie extension} \hbox{ } ~~ \\

Let us begin by a cheap result.

\begin{prop} \label{descNinfty} Corollary \ref{descNar} 1. $(i)$ and 2. $(i)$  holds verbatim if one replaces $Kk_n$ by $K_n$ and $N_{K_{ar}}$ and $N_{K_{ar}}^*$ respectively
by $N_{K_{\infty}}$ and $N_{K_{\infty}}^*$.
\end{prop}
\begin{flushright}$\square$\end{flushright}

To go further we need the following lemma which (partly) is  the derived version of Nakayama's lemma.

\begin{lem}
Let $M\in D^p(\La(G))$, then $M$ is $S$-torsion (ie. $\L_S\Ltens M=0$) if and only $\La(\G)\Ltens_{\La(G)} M\in D^p(\La(\G))$ is finitely generated over $\Zp$.
\end{lem}

Proof. First assume that $M$ is a finitely generated $\La(G)$-module. In this case, the following statements are equivalent:

$(i)$ $M$ is $S$-torsion.

$(ii)$ $M$ is finitely generated over $\La(H)$.

$(iii)$ $\La(\G)\otimes_{\La(G)}M$ is finitely generated over $\Zp$.

$(iv)$ each $Tor_p^{\La(G)}(\La(\G), M)$ is finitely generated over $\Zp$. \smallskip

\noindent

Indeed $(i)\Leftrightarrow (ii)$ is \cite{CFKSV} 2.3,  $(ii)\Leftrightarrow (iii)$ is the topological version of Nakayama's lemma for the ring
$\La(H)$ and the compact $\L(H)$-module $M$ (note that $\La(\G)\otimes_{\La(G)}M\simeq \Zp\otimes_{\La(H)}M$). Finally $(iii)\Leftrightarrow (iv)$ is
an immediate consequence of the fact that $M$ admits a resolution by a perfect complex of $\La(H)$-modules as soon as it is finitely generated over
$\La(H)$.
 \smallskip

Now we come to the general case. As in the case of modules, the implication ``$M$ is $S$-torsion'' $\Rightarrow$ ''$\La(\G)\Ltens_{\La(G)} M$ is
finitely generated over $\Zp$'' is immediate. It thus remains to prove the reciprocal implication. Assume thus that $\La(\G)\Ltens_{\La(G)} M$ is
finitely generated over $\Zp$, ie. that the final terms of the spectral sequence $$Tor_p^{\La(G)}(\La(\G),H^q(M))\Rightarrow Tor_{p-q}(\La(\G),M)$$
are. Then using the above equivalence $(iii)\Leftrightarrow (iv)$ for the module $H^q(M)$ as well as a descending induction on $q$, one gets that
each initial term must also be finitely generated over $\Zp$ (see \cite{Va} 4.1 for details).
\begin{flushright}$\square$\end{flushright}

\begin{prop}\label{S*torsion} Under the assumption ($\mu=0$ - up to iso), $N_{K_\infty}$ and $N_{K_\infty}^*$ are both $S^*$-torsion.
\end{prop}
Proof. Let $A\rightarrow A'$ be an isogeny over $K$ such that $A'$ verifies $(\mu=0)$. Then clearly  $$\Qp\otimes^L_\Zp R\limp R\G^Z(\underline C,T_p(\cal A))\simeq \Qp\otimes^L_\Zp R\limp R\G^Z(\underline C,T_p(\cal A'))$$ and we may thus assume that $A=A'$. By Thm. \ref{La-duality}, $N_{K_\infty}$ (resp. $N_{K_\infty}^*$) is perfect. Taking $R\limp$ along $Kk_\infty/K$ thus gives $$\La(\G)\Ltens_{\La(G)}N_{K_\infty}\simeq N_{K_{ar}} \hbox{\hspace{1cm} (resp.
\hspace{1cm}} \La(\G)\Ltens_{\La(G)}N_{K_\infty}^*\simeq N_{K_{ar}}^*\hbox{\hspace{.5cm})}$$ and the result now follows directly from the previous lemma.
\begin{flushright}$\square$\end{flushright}

\begin{rem} Thanks to the previous proposition, we may form the classes $[N_{K_\infty}]$ and $[N_{K_\infty}^*]$ of  $N_{K_\infty}$ and $N_{K_\infty}^*$  in  $K_0(\mathfrak{M}_G(H))$.
\end{rem}

\begin{cor} Under the assumption ($\mu=0$ - up to iso) each one of the following morphisms of $D^p(\L(G)_{S^*})$
$${\bf 1},\, {\bf 1}-\phi: (I_{K_\infty})_{S^*} \rightarrow  (P_{K_\infty})_{S^*} \hspace{.5cm}, \hspace{1.5cm}
{\bf 1}^*,\, ({\bf 1}-\phi)^*:(P_{K_\infty}^*)_{S^*}\rightarrow (I_{K_\infty}^*)_{S^*}$$ induced by ${\bf 1},\, {\bf 1}-\phi: \underline I\rightarrow
\underline P$, is an isomorphism.
\end{cor}
Proof. For ${\bf 1}$, this is due to the fact that $Lie(D)(-Z)$, and thus $L_{K_\infty}$, is killed by $p$. For ${\bf 1}-\phi$, this is a direct consequence of the previous proposition. \begin{flushright}$\square$\end{flushright}

\begin{rem} It is not the case that $L_{K_\infty}$ is $S$-torsion, because of Prop. \ref{propbc} $(i)$. 
\end{rem}

We are now in position to define the ``$p$-adic $L$-function'' (as well as its dual version) which appears in the main conjecture.

\begin{defn} \label{defl} Under the assumption ($\mu=0$),
the $p$-adic $L$-function  $\cal L=\cal L(A,U,G)$ and it dual version  $\cal L^*=\cal L^*(A,U,G)$ are defined as follows:

$(i)$ $\cal L:=Det(({\bf 1}-\phi){\bf 1}^{-1}|\,(P_{K_\infty})_{S^*})^{-1}\in K_1(\La_{S^*})$.

$(ii)$ $\cal L^*:= Det(({\bf 1}^*)^{-1}({\bf 1}-\phi)^*|\,(P_{K_\infty}^*)_{S^*})^{-1}\in K_1(\La_{S^*})$.
\end{defn}

\begin{rem} \label{remdefL} $(i)$ The notation $\cal L^*$ is justified by the fact that it is deduced from $\cal L$ via the duality involution of $K_1(\La(G)_{S^*})$
(with the convention that duals are viewed as left modules via $g\mapsto g^{-1}$). Indeed $(P_{K_\infty}^*)_{S^*}\simeq
RHom_{\La_{S^*}}((P_{K_\infty}^*)_{S^*},\La_{S^*})$ and $({\bf 1}^*)^{-1}({\bf 1}-\phi)^*$ is $RHom_{\La_{S^*}}(-,\La_{S^*})$-dual to $({\bf 1}-\phi){\bf
1}^{-1}|\,(P_{K_\infty})_{S^*}$.

$(ii)$ Thanks to the previous remark and lemma \ref{rhodual}, we see that for any Artin representation $\rho$ $$\rho(\cal L)=\rho^\vee(\cal L^*)$$
\end{rem}





Let us now review how one classically attaches twisted $L$-functions to $A$ and $U$ via overconvergent isocrystals. By \cite{LST}, the log Dieudonn\'e crystal $D$ over $C^\sharp$  induces an overconvergent $F$-isocrystal $D^\dagger$ over $U/\mathbb F_p$.
Let $L$ be a finite extension of $\Qp$ which is \emph{totally ramified}, $O$ its ring of integers, and consider an Artin representation $\rho:G\to Aut_O(V_\rho)$. Chose $n$ such that $\rho$ factors through $G_n=Gal(U_n/U)$.
 Then we can see $\rho$ as a $O$-representation of
the fundamental group of $U$, having finite (in fact trivial) local monodromy. This representation corresponds to a unique unit-root
overconvergent $F$-isocrystal $U(\rho)^\dagger$ over $U/L$ (see e.g. \cite{Tsz}, 7.2.3) which becomes constant over $U_n$. Consider
$$pr^*:F-iso^\dagger(U/\Q_p)\to F-iso^\dagger(U/L)$$ the natural base change functor from the category of overconvergent $F$ isocrystals over
$U/\Q_p$ to the category of overconvergent $F$ isocrystals over $U/L$. Then we set
$$L(U,A,\rho,s):=L(U,pr^*D^\dagger\otimes U(\rho)^\dagger,p^{-s})$$ where the right hand side
is the classical $L$-function defined in \cite{EL1} associated to the $F$-isocrystal $pr^*D^\dagger\otimes U(\rho)^\dagger$.


\medskip

The rest of this section will be devoted to proving the following theorem, which is our main result:

\begin{thm}\label{main-theorem} The $p$-adic $L$-function and its dual version verify the following properties.

\noindent 1. $(Char)$ In $K_0(\mathfrak{M}_G(H))$, one has the following equalities:

$(i)$ $\partial(\cal L)=[N_{K_\infty}]+[L_{K_\infty}]$.

$(ii)$ $\partial(\cal L^*)=-[N_{K_\infty}^*]-[L_{K_\infty}^*]$.

\noindent 2. $(Interpolation)$: For any totally ramified extension $L$ of $\Qp$ with ring of integers denoted $O$ and every $\rho:G\rightarrow Aut_O(V)$ with contragredient $\rho^\vee:G\rightarrow
Aut_O(V^*)$, one has, in  $L\cup \{\infty\}$:

$(i)$ $\rho(\cal L)=L(U,A,\rho^\vee,1)$.

$(ii)$ $\rho(\cal L^*)=L(U,A,\rho,1)$.
\end{thm}


Proof. 1. $(i)$ According to Lem. \ref{exfrac} there exists a commutative square of the form  $$\diagram{I_{K_\infty}&\hfl{{\bf 1}-\phi}{}&P_{K_\infty}
\cr \vfl{{\bf 1}}{}&&\vfl{}{}\cr
P_{K_\infty}&\hfl{}{}&P_{K_\infty}}$$
whose edges all become isomorphisms after localization by $S^*$. The claimed equality follows from Lem. \ref{frac} since $[Cone({\bf 1}-\phi)]=[N_{K_\infty}[1]]=-[N_{K_\infty}]$ and $[Cone({\bf 1})]=[L_{K_\infty}]$. Here, the
existence of such a diagram is ensured by Lem. \ref{exfrac}.
The proof of $(ii)$ is similar using the dual distinguished triangles describing $N_{K_\infty}^*$ and $L_{K_\infty}^*$.


2. Note that the formula $(ii)$ follows from $(i)$ by Rem. \ref{remdefL} $(ii)$. In order to prove $(i)$ we begin with preliminary results.

\begin{prop}\label{descrig}
There is a canonical isomorphism as follows in the derived category of $L$-vector spaces $$L\Ltens_{\La_O(G)}(
V_{\La_{O}(G)})^*\Ltens_{\La(G)}P_{K_\infty}\simeq R\G_{rig,c}(U/L,pr^*D^\dagger\otimes U(\rho^\vee))$$
(here the $(\L_O(G),\L(G))$-bimodule structure of $(V_{\La_{O}(G)})^*$ is as in the definition of the map $\rho:K_1(\L(G))\to K_1(O)$, see. paragraph \ref{defevrho}). The action of $\phi{\bf
1}^{-1}$ on the left hand term corresponds moreover to the action of Frobenius divided by $p$ on the right hand side.
\end{prop}
Proof. We have isomorphisms as follows:
$$\begin{array}{rcl} L\Ltens_{\La_O(G)} (V_{\La_{O}(G)})^*\Ltens_{{\La_O(G)}}P_{K_\infty}&\mathop{\simeq}\limits^1&
L\Ltens_{\La_O(G)}(V^*\Ltens_\Zp P_{K_\infty})\cr
 &\mathop{\simeq}\limits^2&L\Ltens_{L[G_n]} (V_L^*\Ltens_\Zp P_n)\cr
  &\mathop{\simeq}\limits^3&L\Ltens_{L[G_n]}
(V_L^*\Ltens_L R\G_{rig,c}(U_n/L,pr^*D^\dagger))\cr
&\mathop{\simeq}\limits^4& R\G_{rig,c}(U/L,pr^*D\otimes U^\dagger(\rho^\vee))\end{array}$$ Here the isomorphism $1$ is by Rem. \ref{remtwist} $(i)$, $2$ is obvious and $4$ follows from \'etale cohomological descent for compactly supported rigid cohomology (\cite{CT}) together with the definition of $U(\rho^\vee)$. The isomorphism  $3$ follows from the fact that in   $D(L[G_n])$ we have a Frobenius compatible isomorphism as follows: $$L\Ltens_\Zp R\G(C_n^\s/\Zp,D(-Z))\simeq R\G_{rig,c}(U_n/L,pr^*D^\dagger)$$
Let us explain the latter isomorphism. Since (compactly supported) rigid cohomology is compatible to base change from $\Qp$ to $L$ we may assume that $L=\Qp$ (note that the Frobenius is $L$-linear since $L/\Qp$ is totally ramified) and then apply the following lemma to the diagram having $C^\s_{n}$ as only vertex and edges indexed by $G_n$.

\begin{lem} Consider a diagram (see \cite{TV} Def. 2.1) $X/\Delta$ of proper smooth curves over a perfect field $k$, $Z$ a smooth divisor of $X$, $X^\s=(X,Z)$, $Y=X-Z$, $E$ an $F$-crystal on $X^\s$ and $E^\dagger$ the induced overconvergent $F$-isocrystal (\cite{LST}). If $\Delta$ satisfies the finiteness condition of \cite{TV} Lem. 5.7 then we have a canonical isomorphism $$\Qp\Ltens_\Zp R\G(X^\s/\Zp,E(-Z))\simeq R\G_{rig,c}(Y,E^\dagger)$$ in $D((\Qp\mod)^{\Delta})$. This isomorphism is compatible with Frobenius.
\end{lem}
Proof. By using Poincar\'e duality on both sides of the isomorphism, we are reduced to establish an isomorphism $$\Qp\Ltens_\Zp R\G(X^\s/\Zp,E)\simeq R\G_{rig}(Y,E^\dagger)$$
in $D((\Qp\mod)^{\Delta^{op}})$. Moreover the compatibility with Frobenius is a formal consequence of the rest since one can replace the diagram $X$ by the diagram $Frob:X\to X$ of type $\Delta\times [1]$.
Using cohomological descent on both the crystalline and the rigid side (see \cite{CT}), we can always assume given an exact immersion $X^\s\subset P^\s$ into a log smooth formal log scheme $P^\s/W(k)$ as in \cite{TV} Lem. 5.7. Then, after tensorisation by $\Q_p$, the crystalline cohomology of $E$ can be computed (see \cite{Sh} and \cite{Tr}, proof of proposition 3.3) as the log de Rham cohomology on the tube $]X[_P$ of the associated module with connection $E_K$ ($K=Frac(W(k)$) and maps naturally to the de Rham cohomology on $V$, some strict neighborhood of $]Y[_P$ in $]X[_P$, of $E_K|_V$, which is nothing but the rigid cohomology of $E^\dagger$. We have then constructed a map $$\Qp\Ltens_\Zp R\G(X^\s/\Zp,E)\to R\G_{rig}(Y,E^\dagger).$$
To show that this map is an isomorphism, we need to show that this is the case on each vertex of the diagram.  Now we can always assume given a proper smooth lifting of each vertex (they need not be compatible) and we conclude using \cite{LST}, 4.2 and \cite{Tsj}, 1.5.
\begin{flushright}$\square$\end{flushright}


\begin{prop} \label{weil} \cite{EL2}
Let $E^\dagger\in F-iso^\dagger(U/\Q_p)$.
$$L(U,E^\dagger,\rho,t)=\prod_{i=0}^2det_L(1-tF|H^i_{rig,c}(U/L,pr^*D^\dagger\otimes U(\rho)^\dagger)^{(-1)^{i+1}}$$
\end{prop}
\begin{flushright}$\square$\end{flushright}

Now we can prove the interpolation formula 2.$(i)$ of Thm. \ref{main-theorem}. For simplicity, we set
$${_RP}:=R\Ltens_{\La_O(G)}(V_{\La_O(G)})^*\Ltens_{\La(G)}P_{K_\infty}\hbox{\hspace{.5cm}and\hspace{.5cm}}
{_RH^q}:=H^q(_RP)$$ for any $\La_O(G)$-algebra $R$.
Note that ${_{R'}}H^q\simeq R'\otimes_R{_R}H^q$ if $R'/R$ is flat.

Applied to $P_{K_\infty}$, Lem. \ref{rho_Det} gives \begin{eqnarray}\label{pfmn1}\rho_{Q_O(\G)}\cal L=\prod_qdet_{Q_O(\G)}(({\bf 1}-\phi){\bf 1}^{-1}|
_{Q_O(\G)}H^q)^{(-1)^{q+1}}\end{eqnarray} and Lem. \ref{Tor} applied to ${_{\La_O(\G)_I}H^q}$ gives:
\begin{eqnarray}\label{pfmn2}\begin{array}{ccrcl}&&\epsilon det_{Q_O(\G)}(({\bf
1}-\phi){\bf 1}^{-1}|{_{Q_O(\G)}H^q})&=&
det_L(({\bf 1}-\phi){\bf 1}^{-1}|Tor_0^{\La_O(\G)_I}(L,{_{\La_O(\G)_I}H^q}))\\ &&&&.det_L(({\bf 1}-\phi){\bf
1}^{-1}|Tor_1^{\La_O(\G)_I}(L,{_{\La_O(\G)_I}H^q}))^{-1}\end{array}\end{eqnarray} whenever the right hand side makes sense (recall that the map $\epsilon$ was defined in (\ref{defrho})).

Now Prop. \ref{descrig} shows that $$\begin{array}{rcl} R\G_{rig,c}(U/L,pr^*D^\dagger\otimes U(\rho^\vee))&\simeq& {_LP}\\
&\simeq& L\Ltens_{\La_O(\G)_I} {_{\La_O(\G)_I}P}\end{array}$$ Whence a spectral sequence which degenerates into short exact sequences

\begin{eqnarray}{\small\label{pfmn3}\diagram{&&&&&H^{2}_{rig,c}(U/L,pr^*D^\dagger\otimes U(\rho^\vee))&\surjfl{}{}&
Tor_0^{\La_O(\G)_I}(L,{_{\La_O(\G)_I}H^2})\cr &&&Tor_1^{\La_O(\G)_I}(L,{_{\La_O(\G)_I}H^2})&\injfl{}{}&H^{1}_{rig,c}(U/L,pr^*D^\dagger\otimes
U(\rho^\vee))&\surjfl{}{}& Tor_0^{\La_O(\G)_I}(L,{_{\La_O(\G)_I}H^1})\cr
&&&Tor_1^{\La_O(\G)_I}(L,{_{\La_O(\G)_I}H^1})&\injfl{}{}&H^{0}_{rig,c}(U/L,pr^*D^\dagger\otimes U(\rho^\vee))&\surjfl{}{}&Tor_0^{\La_O(\G)_I}(L,{_{\La_O(\G)_I}H^0})\cr
&&&Tor_1^{\La_O(\G)_I}(L,{_{\La_O(\G)_I}H^0})&=&0& }}
\end{eqnarray}


{\bf Claim}. The operator $({\bf 1}-\phi){\bf 1}^{-1}$ is bijective on each $Tor_{-i}^{\La_O(\G)_I}(L,{_{\La_O(\G)_I}H^j}))$ except maybe for
$(i,j)=(0,1)$. Moreover, the operator $1-p^{-1}F$ is bijective on $H^{q}_{rig,c}(U/L,pr^*D^\dagger\otimes U(\rho^\vee))$ except maybe for $q=1$. \smallskip

Let us prove this claim. First we observe that the second statement implies the first one since by Nakayama's lemma, $({\bf 1}-\phi){\bf 1}^{-1}$ is bijective on ${_{\La_O(\G)_I}}H^j$ if and
only if it is on $Tor_{0}^{\La_O(\G)_I}(L,{_{\La_O(\G)_I}H^j})$, in which case it is also bijective on
$Tor_{1}^{\La_O(\G)_I}(L,{_{\La_O(\G)_I}H^j}))$. Now the second claim follows from the long exact sequence (cf. Thm. \ref{key-dt} and the isomorphism between the third and last term at the beginning of the proof of Prop. \ref{descrig})
$$\diagram{L\otimes_{L[G_n]} (V_L^*\Ltens_\Zp N_n)&\hflcourte{}{}& H^{q}_{rig,c}(U/L,pr^*D^\dagger\otimes U(\rho^\vee))&\hflcourte{1-p^{-1}F}{}& H^{q}_{rig,c}(U/L,pr^*D^\dagger\otimes U(\rho^\vee))&\hflcourte{}{}&\dots}$$
since by Cor. \ref{arcoh} $\Q\otimes N_n$ is concentrated in degrees $[1,2]$.
\smallskip

We may now conclude since the claim ensures that no indeterminate product of the form $0.\infty$ occurs in:

- the formula obtained by combining (\ref{pfmn1}) with (\ref{pfmn2}) (note that the left hand term is nothing but $\rho(\cal L)$.

{\small
$$\label{pfmn4}\epsilon \rho_{Q_O(\G)}(\cal L)\prod_qdet_{Q_O(\G)}(({\bf 1}-\phi){\bf 1}^{-1}|
_{Q_O(\G)}H^q)^{(-1)^{q+1}}=
\prod_{i,j}det_L(({\bf 1}-\phi){\bf 1}^{-1}|Tor_{-i}^{\La_O(\G)_I}(L,{_{\La_O(\G)_I}H^j}))^{(-1)^{i+j+1}}$$}

- the formula given by the exact sequences (\ref{pfmn3})
{\small $$\prod_{i,j}det_L(({\bf 1}-\phi){\bf 1}^{-1}|Tor_{-i}^{\La_O(\G)_I}(L,{_{\La_O(\G)_I}H^j}))^{(-1)^{i+j+1}}=
\prod_qdet_L(1-p^{-1}F|H^{q}_{rig,c}(U/L,pr^*D^\dagger\otimes U(\rho^\vee))^{(-1)^{(q+1)}}$$}

- the formula given by Lem. \ref{weil}
{\small $$\prod_qdet_L(1-p^{-1}F|H^{j}_{rig,c}(U/L,pr^*D^\dagger\otimes U(\rho^\vee))^{(-1)^{(q+1)}}= L(U,A,\rho^\vee,1)\hbox{.}$$}
\begin{flushright}$\square$\end{flushright}


%
%

\begin{rem}
$(i)$ In the above proof we have retrieved the well known fact that $\rho(\cal L)=L(U,A,\rho^\vee,1)$ may take the value $0$, but not $\infty$.

$(ii)$ In \cite{LLTT}, the author considers the complex $RHom_\Zp(\underline L,\Qp/\Zp)=\underline L^*[1]$ instead of $\underline L^*$. This
explains the difference of signs between formula 1. $(ii)$ and those  of loc. cit.
\end{rem}

%
%


\end{document}